\documentclass[final,1p,times]{elsarticle}

\usepackage{amssymb}
\usepackage{bm}
\usepackage{amsthm}
\usepackage{amscd}
\usepackage{amsmath}
\usepackage{amsfonts}
\usepackage{amssymb}
\usepackage{graphicx}
\usepackage{color}
\newtheorem{theorem}{Theorem}

\newtheorem{lemma}[theorem]{Lemma}

\newtheorem{definition}[theorem]{Definition}
\newtheorem{remark}[theorem]{Remark}
\usepackage{mathrsfs}
\usepackage{titletoc}

\usepackage{arydshln}
\newcommand{\be}{\begin{equation}}
	
	\newcommand{\ee}{\end{equation}}
\newcommand{\bea}{\begin{eqnarray}}
	\newcommand{\eea}{\end{eqnarray}}
\newcommand{\bna}{\begin{eqnarray*}}
	\newcommand{\ena}{\end{eqnarray*}}

\renewcommand{\le}{\left}
\newcommand{\ri}{\right}

\journal{***}

\begin{document}
	
	\begin{frontmatter}
		\title{Free boundary regularity in  nonlinear one-phase Stefan problem 
			}

		\author[ruc1]{Yamin Wang}
		\ead{ym.wang@nus.edu.sg}
		\address[ruc1]{Department of Mathematics, National University of Singapore, 119076, Singapore.}
		%\cortext[cr]{Corresponding author.}
		%\tnotetext[sw]{The author is supported by the China Scholarship Council in 2021  and the Outstanding Innovative Talents Cultivation Funded Programs 2021 of Renmin University of China.}		
		%\cortext[cr]{Corresponding author.}

		\begin{abstract}
	We study  the regularity of the free boundary in  one-phase Stefan problem with nonlinear operator. Using the Hodograph transform and a  linearization technique, we prove  that flat  free boundaries are $C^{1, \alpha}$ in space and time. When  the  operator is  concave (or convex) and smooth, the free boundary  is smooth. %This is an extension  of a regularity result derived by De Silva-Forcillo-Savin who treated the linear parabolic counterpart. 
		\end{abstract}
		
		\begin{keyword} 
			Nonlinear parabolic equations,  Stefan problem,  free boundary, regularity. 		%\MSC[2010] %58J05\sep 58J32
			
		\end{keyword}
		
	\end{frontmatter}
	
	\section{Introduction}
	We consider in this paper the one-phase  Stefan problem. More precisely,  let $\Omega\subset \mathbb{R}^n$ be a bounded domain, and let $u=u(x,t)$ denote the temperature of the medium at a point $x\in \Omega$ at $t\in \mathbb{R}^+:=[0,+\infty)$.   The \emph{classical Stefan problem} can be  formulated as follows:
\be\label{ak1}\le\{\begin{array}{lll}
	u_t=\Delta u\quad\;\; \; \text{in}\quad\;\Omega^+_u,\\[1ex]
	u_t=|\nabla u|^2\quad \text{on}\quad  \partial\Omega^+_u,\\[1ex]
	u\geq 0\quad \quad \;\;\; \text{in}\quad\; \Omega\times [0,T],
\end{array}\ri.\ee
where  the positive set 
$$\Omega^+_u:=\{(x, t)\in \Omega\times (0, T]: u(x,t)>0\},$$
and the free boundary
 $$ \partial\Omega^+_u:=\{(x,t)\in \Omega \times(0, T]: u(x,t)=0\}.$$ 
Moreover,  the Laplacian and the gradient with  regard to the spatial variables are denoted by  $\Delta, \nabla $ respectively. It is well known that the classical  Stefan problem (\ref{ak1}) describes the phase transition between solids and liquids, typically the melting of ice to water, for example \cite{fri,me,ru}. In this model,  $\Omega^+_u$ represents the water, while its complement  describes the region of unmelted  ice.  \vspace{0.2cm} 
 
One notices that the second  condition in (\ref{ak1}) determines the velocity of the moving interphase. That is, if $u$ is smooth up to the boundary, then the free boundary moves with the normal velocity $V=u_t/|\nabla u|$. Therefore, the second relation indicates 
$$V=\nabla u\cdot \nu=|\nabla u|\quad \text{on}\quad \partial\Omega^+_u, $$
where $\nu=\nabla u/|\nabla u|$ denotes the spatial unit normal vector of $\partial\Omega^+_u$ at $(x, t)$. \vspace{0.2cm}

In the above Stefan problem, the regularity of  solutions and the free boundaries is of particular  interest.  For what concerns the solution,  one can establish the optimal regularity plus the quantitive properties such as nondegeneracy and semi-convexity in space (see \cite{cc2,cc1,cc3}). However, a central mathematical challenge is to understand the geometry or regularity of the free boundary. It also plays an important role in proving further regularity of the solutions.  This theory was developed by Caffarelli in the groundbreaking paper \cite{cc2}, in which he proved that the  free boundary is smooth outside of a certain set of cusp-like singularities. It was later derived by Figalli-Ros-Oston-Serra \cite{figar} that the singular set has parabolic Hausdorff dimension at most $n-1$.  \vspace{0.2cm}

There has been  an extensive literature on the regularity of the free boundary for the classical Stefan problem (\ref{ak1}).  In the celebrated work,  Athanasopoulos-Caffarelli-Salsa \cite{acs1} showed that Lipschitz free boundaries of the two-phase Stefan problem are smooth under  a nondegeneracy condition. After that in  \cite{acs2}, they obtained the same results for flat free boundaries.  The general idea of the contribution \cite{acs1, acs2}  was inspired by the elliptic case  \cite{ca1,ca2}. Using a different method estabished in the elliptic counterpart \cite{desilva}, De Silva-Forcillo-Savin  \cite{r1} recently established  an equivalent  conclusion to the  flatness result  of  \cite{acs2}. \vspace{0.2cm}

As counter example in \cite{acs1} (see also \cite{ss1,cas1} for the one-phase case), Lipschitz free boundary in evolution problems does not enjoy instantaneous regularization. In general, a Lipschitz free boundary could exhibit a hyperbolic behavior, i.e,  a corner can persist for an amount of time.  Thus a nondegeneracy condition would be indispensable in dealing  with this problem.  Concerning this phenomenon, Choi-Kim \cite{sck1}  proved that the one-phase Stefan free boundary regularizes in space and time when starting from Lipschitz initial free boundary with small Lipschitz constant.   \vspace{0.2cm} 

 In this work,  we study the problem (\ref{ak1}) with  fully nonlinear operator, namely the following  \emph{nonlinear Stefan problem} 
\be\label{a1}\le\{\begin{array}{lll}
	u_t=\mathcal {F}(D^2 u)\quad\;\, \text{in}\quad\Omega^+_u,\\[1ex]
	u_t=\mathcal{G}(|\nabla u|)\quad\; \;\,\text{on}\quad \partial\Omega^+_u,\\[1ex]
		u\geq 0\quad \quad\quad \quad \;\,\text{in}\quad\; \Omega\times [0,T],
\end{array}\ri.\ee
where  $D^2u$ stands for the spacial Hessian of $u$ and $\mathcal {F}$ is the fully nonlinear elliptic operator. Additionally, we assume that $\mathcal {F}$ and  $\mathcal{G}$ satisfy the  conditions
 \be\label{addf1}\le\{\begin{array}{lll}
   	\mathcal {F}\;\text{is\;uniformly\;elliptic},\;\mathcal {F}(0)=0,\\[1ex]
   \partial_p (\mathcal{G}(p)/ p)\geq c>0,\; \mathcal{G}\geq 0,
 \end{array}\ri.\ee
where $c$ is a constant. In our generality, the regularity of free boundaries is an object of investigations. We shall focus on the  perturbative estimates for this nonlinear  problem (\ref{a1}).\vspace{0.2cm}

The free boundary regularity with nonlinear operator $\mathcal{F}$ has attracted much attention in the last decades. In elliptic case, when $\mathcal{F}$ is homogeneous of degree one,  several authors extended the results of seminal works  \cite{ca1,ca2} to various kinds of nonlinear
operators.  For works in this direction,  Wang \cite{w1,w2} considered the concave operator  of
the form $\mathcal{F}=\mathcal{F}(D^2 u)$ and then it was generalized  by Feldman \cite{feld} to a class of operators $\mathcal{F}=\mathcal{F}(D^2 u, Du)$ by removing concavity assumption.  Without the concavity or homogeneity assumption of $\mathcal{F}$, De Silva-Ferrari-Salsa \cite{ss1} established the $C^{1, \alpha}$ regularity for flat free boundary in the problems with distributed source. We also refer to \cite {af,fe1} for the type with H$\ddot{\text{o}}$lder dependence on $x$ of the operator.  In parabolic case,  Milakis \cite{Milakis, Milakis2} considered a two-phase problem  with $\mathcal{F}$ concave and homogeneous of degree 1. \vspace{0.2cm}

We shall discuss about the regularity of the  free boundaries of the nonlinear problem (\ref{a1}) under  a flatness hypothesis. The ideas for this are closer in spirit to the pioneers' work \cite{r1} concerning the classical setting (\ref{ak1}).  Parallel with \cite[Theorem 1.2]{r1}, our main result roughly states that a sufficiently flat solution $u$ to (\ref{a1}) in a certain ball in space and time has  $C^{1,\alpha}$   free boundary in the interior.  To clarify the notation of flatness that is often used, we  let $\Omega=B_{r}(x)$ and for reasons of convenience,
$$\mathcal{Q}_{\lambda}:=B_{\lambda}\times [-\lambda, 0],$$
where $B_r(x)$ stands for the ball in $\mathbb{R}^n$ of radius $r$ centered at $x$ (the dependence on $x$ will often be omitted if $x=0$). Let
\begin{equation}\label{a3}
	l_{a,b}(x,t):=a(t)\cdot x+b(t), \quad x\in \mathbb{R}^n,
\end{equation}
where $n$-dimensional space variable $x=(x_1,x_2,\ldots, x_n)$ and $a(t)=(a_1, \ldots, a_{n-1}, a_n(t))$ with $a_i\in \mathbb{R}$ for $ i=1, \ldots, n-1.$  We say that  $u(x,t)$ is \emph{$\epsilon_0$-flat} in a ball of size $\lambda$ in space and time, if  $u$ is trapped in a strip of width $\epsilon_0 \lambda$ by two parallel hyperplanes, i.e. 
\begin{equation}\label{a2}
	(a_n(t)x_n-b(t)-\epsilon_0 \lambda)^+\leq u(x,t)\leq (a_n(t)x_n-b(t)+\epsilon_0\lambda)^+\quad \text{in}\quad \mathcal{Q}_{\lambda}.
\end{equation}

Our  rigorous statement  is below. The notion of viscosity solution is reviewed in (Section 2, Definition \ref{de3}).  
\begin{theorem}\label{t1}
	For $\mathcal {F}, \mathcal{G}$ satisfying  (\ref{addf1}), 	let $u$ be a viscosity solution to (\ref{a1})   in $\mathcal{Q}_{\lambda}$ for some $\lambda\leq 1$.  Assume  that $(0,0)\in \partial\Omega^+_u$ and for fixed constant $M>1$, there are universal small constants $\bar{\epsilon}, c_0>0$ such that $u$ is $\bar{\epsilon}$-flat with $b'(t)=-\mathcal{G} (a_n)$ and 
	\begin{equation}\label{ar2}
		M^{-1}\leq a_n\leq M,\quad \quad |a^{\prime}_n|\leq c_0\lambda^{-2}.
	\end{equation}
	Then in $ \mathcal{Q}_{\lambda/2}$ the free boundary $\partial\Omega^+_u$ is a $C^{1,\alpha}$ graph in the  direction of $x_n$. If in addition that $\mathcal {F}$  is convex or concave and smooth, then $\partial\Omega^+_u $ is smooth in $\mathcal{Q}_{\lambda/2}$.
\end{theorem}

Here and henceforth, we say  that a constant is \emph{universal} if it depends only on the dimension $n, M$ and ellipticity constants. It is worth to mention that our assumptions in the above theorem  implies  the nondegeneracy property of $u$  (see Section 5 for more details).   \vspace{0.2cm}

Now let us make some initial remarks on the  proof Theorem \ref{t1}.  The main strategy in the proof  is  to show that the graph of $u$ enjoys an improvement of flatness property. Then the $C^{1, \alpha}$ regularity in space and time of the free boundary is achieved by an iteration procedure.  For this purpose,  the first step we shall proceed relies on the hodograph transform, which passes  (\ref{a1})   into an equivalent problem with fixed  boundary (see (\ref{a9}) below). This is due to the  lack of natural rescaling for the equation. To illustrate the idea, we take $\mathcal {F}(D^2 u)=\Delta u$ and  $\mathcal {G}(|\nabla u|)=|\nabla u|^2$ for example. Let the function $u$ solve (\ref{a1})  in $\mathcal{Q}_{\lambda}$. To preserve the equation in positive phase, one uses the  parabolic rescaling
$$u_{\lambda}(x,t)=\frac{u(\lambda x, \lambda^2 t)}{\lambda}, \quad (x,t)\in B_{1}\times [-\lambda^{-1}, 0].$$ 
 The  letting $\lambda\rightarrow 0$ formally, we find that  $w:=\lim_{\lambda\rightarrow0}u_{\lambda}$ solves 
 $$
	w_t=0\quad  \text{on}\quad (B_{1}\times (-1, 0])\cap \{w=0\}.
$$
Thus the free boundary condition degenerates. If we  adopt  the hyperbolic rescaling  $$u_{\lambda}(x,t)=\frac{u(\lambda x, \lambda t)}{\lambda},\quad (x,t)\in \mathcal{Q}_{1},$$ 
as $\lambda$ goes to $0$, the limiting solution sloves 
$$
w_t=0\quad  \text{in}\quad  (B_{1}\times (-1, 0])\cap \{w>0\},
$$
which is the so-called Hele-Shaw problem (see e.g. \cite{cg,ckim}). However the good continuity in time for the solution is failed.  \vspace{0.2cm}

Secondly,  the improvement of flatness is obtained via a suitable compactness and linearization argument. More concretely, the  nonlinear equation is linearized to an oblique derivative parabolic problem, for which various regularity estimates can be  proved.  In this process, the key ingredient is to establish a diminishing of oscillation property.  A useful tool  in proving this will be Harnack type inequalities for the solutions to general equations with the same type of measurable coefficients.  Although we borrow perturbation techniques  from  \cite{r1}, there are delicate difficulties and challenges  arising from the fully nonlinear term to overcome.  \vspace{0.2cm}

Before ending this introduction, we would like to mention  \cite{da} for higher regularity of the free boundary in  nonlinear Stefan problem, and \cite{ak,fgalli} for local structure of the free boundary in parabolic obstacle  problem.   \vspace{0.2cm}

The remaining part of this paper is structured as follows. In Section 2, we provide notations and definitions used throughout the paper and  present 
auxiliary results which will be applied in the proof of Theorem \ref{t1}.  In Section 3, we  perform the Hodograph transform and perturbative arguments for both the linear and nonlinear problems.  In Section 4, we  establish an improvement flatness result,  while in Section 5, we complete the proof of  Theorem \ref{t1} by applying this result and Schauder estimates.

\section{Preliminaries}		
In this section we recall the notion of viscosity solutions and present some known results about nonlinear parabolic equations. Firstly, we talk about the fully nonlinear elliptic operator  and refer to \cite{cx,cs,gl} for a comprehensive treatise and introduction.  Let $\mathcal{S}$ denote the space of the real $n\times n$ symmetric matrices.  Assume that $\mathcal {F}: \mathcal{S}\rightarrow \mathbb{R}$ is \emph{uniformly elliptic}, i.e., there exists a  constant $1\leq \Lambda<+\infty$ such that for any $\mathcal{M},\mathcal{N}\in \mathcal{S}$ with positive semi-definite $\mathcal{N}$,
\begin{equation}\label{add1}
	\Lambda^{-1} \|\mathcal{N}\|\leq \mathcal {F} (\mathcal{M}+\mathcal{N})-\mathcal {F} (\mathcal{M})\leq \Lambda \|\mathcal{N}\|,
\end{equation}
where  $\|\mathcal{M}\|$  denotes the $(L^2, L^2)$-norm of $\mathcal{M}$, i.e. $\|\mathcal{M}\|=\sup_{|x|=1}\|\mathcal{M} x\|$.  And assume $\mathcal {F}(0)=0$, which  is not essential since we can consider $\hat{\mathcal {F}}(D^2 u):=\mathcal {F}(D^2 u)-\mathcal {F}(0)$ which fulfills (\ref{add1}) with the same ellipticity constants.  \vspace{0.2cm}

In general, solutions of (\ref{a1}) develop singularities in time, so classical solutions may not be expected to exist globally in time. We employ the notation of viscosity solutions introduced in  \cite{cl}. 
To this end, we need the following standard notion.

\begin{definition}\label{de1} (Contact)
We say that a function $\phi$ touches a function $u$ by above (resp. below) at $(x_0,t_0)$ in a parabolic cylinder $E_r(x_0, t_0):=B_r(x_0)\times (t_0-r^2,t_0]$, if 
$\phi(x_0,t_0)=u(x_0,t_0)$ and 
$$u(x,t)\leq \phi(x,t)\quad  (\text{resp.}\; u(x,t)\geq \phi(x,t))$$
 for all $(x,t)\in E_r(x_0, t_0)$. If the above inequality is strict in $ E_r(x_0, t_0)\backslash \{x_0, t_0\}$, we say that $\phi$ touches $u$ strictly from above (resp. below).
\end{definition}

Formally, viscosity solutions are the functions that satisfy a local comparison principle on parabolic neighborhoods with barriers which are the classical solutions of the problem.    We denote that $u\in C^{2,0}_x\cap C^{0,1}_t$ if $u$ is $C^2$ continuous with respect to $x$ and is $C^1$ continuous with respect to $t$. 

\begin{definition}\label{de2} (Comparison solution)
A continuous function $ \phi$ is a comparison supersolution (resp. subsolution) of (\ref{a1})  if $\phi\in  C^{2,0}_x\cap C^{0,1}_t$,  $\mathcal{G}(|\nabla \phi|)\not=0$ and 
\begin{equation*}
	\le\{\begin{array}{lll}
	\phi_t\geq\mathcal {F}(D^2 \phi)\quad \text{in}\quad\Omega^+_\phi,\\[1.5ex]
	\phi_t\geq\mathcal{G}(|\nabla \phi|)\quad \text{on}\quad \partial\Omega^+_\phi.
\end{array}\ri.
	\end{equation*}
If the above inequalities are strict, we say that $ \phi$ is a strict supersolution (resp. subsolution). We say that $ \phi$ is a classical solution to a free boundary problem if it is both subsolution and supersolution.
\end{definition}

\begin{definition}\label{de3} (Viscosity solution) An upper (resp. lower) semi-continuous function $u\geq 0$ in  $\Omega \times [0, T]$ is called a viscosity subsolution (resp. supersolution) of  (\ref{a1})  if the following conditions are satisfied: \vspace{0.2cm} \\
(i) $u_t\geq \mathcal {F}(D^2 u)$ (resp. $u_t\leq \mathcal {F}(D^2 u)$) in $\Omega^+_u$ in the viscosity sense, i.e. if $\phi\in C^{2,0}_x (\Omega^+_u) \cap C^{0,1}_t(\Omega^+_u)$ touches $u$ from below  (resp. above) at $(x_0,t_0)\in \Omega^+_u$, then 
$$\phi_t (x_0,t_0)\geq\mathcal {F}(D^2 \phi(x_0,t_0))\quad (\text{resp.}\; \phi_t (x_0,t_0)\leq\mathcal {F}(D^2 \phi(x_0,t_0))).$$
(ii) Any strict comparison supersolution $\phi^+$ (resp. subsolution) cannot touch $u$ by above (resp. below) at any point $(x_0,t_0)\in \partial\Omega^+_\phi$.   \vspace{0.2cm}

A function $u\in C(\Omega \times [0, T])$ is called a viscosity solution of (\ref{a1}), if it is both a viscosity subsolution and a viscosity supersolution.
\end{definition}
\noindent \begin{remark}\label{rmk2} \emph{ As usually, the above definition is equivalent  to the one in which $\phi\in C^{2,0}_x  \cap C^{0,1}_t$ is replaced by  a paraboloid 
\begin{equation*}\label{vtrr22}
	P(x,t)=\frac{1}{2}x^TAx+Bt+Cx+D
\end{equation*}
with $A, B, C, D$ are some constants. }\end{remark}

For any symmetric matrix $\mathcal{R}\in \mathcal{S}$ and parameters $\Lambda, \Lambda^{-1}$, we consider the extremal Pucci operators $\mathcal{P}^{\pm}:\mathbb{R}^{n\times n}\rightarrow\mathbb{R}$ (see e.g.  \cite{cx,gl}) 
$$\mathcal{P}^{-}\Big(\mathcal{R}, \Lambda^{-1},\Lambda\Big):=\Lambda^{-1}\sum_{e_i>0}e_i+\Lambda\sum_{e_i<0}e_i,\quad\quad \mathcal{P}^{+}\Big(\mathcal{R}, \Lambda^{-1},\Lambda\Big):=\Lambda^{-1}\sum_{e_i<0}e_i+\Lambda\sum_{e_i>0}e_i,$$
where $e_i=e_i(\mathcal{R})$ are the eigenvalues of $\mathcal{R}$. It is easy to see that
\begin{equation}\label{p1}
	\mathcal{P}^{-}\Big(\mathcal{R}, \Lambda^{-1},\Lambda\Big)=\inf_{P\in \mathfrak{P}_{ \Lambda^{-1},\,\Lambda}} \{\text{tr}(P \mathcal{R})\},\quad 
	\mathcal{P}^{+}\Big(\mathcal{R}, \Lambda^{-1}, \Lambda\Big)=\sup_{P\in \mathfrak{P}_{ \Lambda^{-1},\,\Lambda}} \{\text{tr}(P \mathcal{R})\},
\end{equation}
where $\mathfrak{P}_{ \Lambda^{-1},\,\Lambda}:=\{P\in  \mathcal{S}\,|\,\Lambda^{-1} \text{I}\leq P\leq \Lambda\text{I}\}$. We know that the Pucci operators enjoy some basic properties as follows. 
\begin{theorem}\label{thf1} (Properties of Pucci)
	It holds that\\\vspace{0.1cm} 
(i) For $1\leq \Lambda_1\leq \Lambda_2$, 
\begin{equation*}\label{pe1}
	\mathcal{P}^{-}\Big(\mathcal{R}, \Lambda_2^{-1}, \Lambda_2\Big)\leq 	\mathcal{P}^{-}\Big(\mathcal{R}, \Lambda_1^{-1}, \Lambda_1\Big), \quad 	\mathcal{P}^{+}\Big(\mathcal{R}, \Lambda_1^{-1}, \Lambda_1\Big)\leq 	\mathcal{P}^{+}\Big(\mathcal{R}, \Lambda_2^{-1}, \Lambda_2\Big).
\end{equation*}
(ii) For $\mathcal{R}_1, \mathcal{R}_2\in \mathcal{S}$,  
$$\mathcal{P}^{-}\Big(\mathcal{R}_1, \Lambda^{-1}, \Lambda\Big)+ \mathcal{P}^{-}\Big(\mathcal{R}_2, \Lambda^{-1}, \Lambda\Big) \leq \mathcal{P}^{-}\Big(\mathcal{R}_1+ \mathcal{R}_2, \Lambda^{-1}, \Lambda\Big), $$
and 
$$ \mathcal{P}^{+}\Big(\mathcal{R}_1+ \mathcal{R}_2, \Lambda^{-1}, \Lambda\Big) \leq \mathcal{P}^{+}\Big(\mathcal{R}_1, \Lambda^{-1}, \Lambda\Big)+ \mathcal{P}^{+}\Big(\mathcal{R}_2, \Lambda^{-1}, \Lambda\Big).$$
(iii) If $\alpha\geq 0$
$$ \mathcal{P}^{\pm}\Big(\alpha\mathcal{R}, \Lambda^{-1}, \Lambda\Big)=\alpha \mathcal{P}^{\pm}\Big(\mathcal{R}, \Lambda^{-1}, \Lambda\Big).$$ 

\end{theorem}

Next we present the Harnack inequality  \cite[Theorem 4.32]{cs} for solutions to fully nonlinear equations.  Some notations are set up as follows
\begin{equation}\label{frq1}
	\begin{split}
	&	\mathcal{C}_r(x_0, t_0):=B_r(x_0)\times (t_0-r^2,t_0+r^2),\quad 	\mathcal{C}^-_{r/2}(x_0, t_0):=B_{r/2}(x_0)\times \Big(t_0-r^2,t_0-\frac{r^2}{2}\Big),\\
	&	\mathcal{C}^+_{r/2}(x_0, t_0):=B_{r/2}(x_0)\times \Big(t_0+\frac{r^2}{2},t_0+r^2\Big),\quad 	\mathcal{C}_r^{\lambda}:=B_{r}(x_0)\times (t_0- \lambda r^2, t_0+ \lambda r^2).
	\end{split}
\end{equation}

\begin{theorem}\label{th1}(Harnack inequality)
For $r>0$ and  $(x_0,t_0)\in  \mathbb{R}^{n+1}$,  let $u\geq 0$ satisfy
\be\label{ah9}\le\{\begin{array}{lll}
	u_t\leq \mathcal{P}^{+}\left(D^2 u, \Lambda^{-1}, \Lambda\right)+\Upsilon_0\quad \text{in}\quad \mathcal{C}_r(x_0, t_0),\vspace{0.2cm} \\
	u_t\geq \mathcal{P}^{-}\left(D^2 u, \Lambda^{-1}, \Lambda \right)-\Upsilon_0\quad \text{in}\quad \mathcal{C}_r(x_0, t_0),
\end{array}\ri.\ee
where $\Upsilon_0\geq0$ is some constant. Then there is constant $\tilde{C}>0$ depending only on $n, \Lambda$ and  $\Lambda ^{-1}$ such that 
\begin{equation*}
	\mathop{\sup\;\;u}_{\mathcal{C}^-_{r/2}(x_0, t_0)}\leq \tilde{C} \Big(\mathop{\inf\;\; u}_{\mathcal{C}^+_{r/2}(x_0, t_0)}+r^2 \Upsilon_0\Big).
		\end{equation*}
\end{theorem}

The next theorem  states that the infimum of a supersolution $u$ can be controlled  by the $L^p$ norm, where $\mathop{\fint }_{Q}:=\frac{1}{|Q|}\int_{Q}$. For its proof we can consult  \cite[Theorem 4.15]{cs}  and  \cite[Corollary 4.14]{wlh1} .  
\begin{theorem}\label{thw1}(Weak harnack inequality)
For $r>0$ and  $(x_0, t_0)\in \mathbb{R}^{n+1}$, let $u>0$ satisfy
	\begin{equation*}
		u_t\geq \mathcal{P}^{-}\Big(D^2 u, \Lambda^{-1}, \Lambda\Big) \quad \text{in}\quad \mathcal{C}_r(x_0, t_0).
	\end{equation*}
Then there is universal constant $C>0$ and $p\in (0,1)$ such that
$$   \Big(	\mathop{\fint \;\;u^p}_{\mathcal{C}^-_{r/2}(x_0, t_0)}\Big)^{1/p}\leq \tilde{C} \mathop{\inf \;\; u}_{\mathcal{C}^+_{r/2}(x_0, t_0)}.$$
\end{theorem}

A regularity property for the solutions to linear parabolic equation that is coming  is due to  \cite[Section 2.4]{cs}. 
\begin{theorem}\label{th3}(Regularity estimates)
	Let $u$ be solution of $u_t=\Delta u$ in $B_R\times(0, T)$.  Then
	there holds for any multi-index $|\alpha|=k$, 
	$$\left|D^\alpha u(0)\right|\leq \frac{C}{R^k}\max_{B_R\times(0, T)}|u|.$$
	Here $C$ is a constant that depends on $n$ and $k$.
\end{theorem}

 Consider the parabolic equation with fully nonlinear operator
\begin{equation}\label{vtr22}
	u_t=\mathcal {F}(D^2 u)\quad \text{in}\quad B_R\times (0, T).
\end{equation}
The following comparision principle related to (\ref{vtr22}) appears in \cite[Theorem 14.1]{gl}.

\begin{theorem}\label{th2}(Comparision principle)
Assume that $u$ and $v$ are respectively  subsolution and supersolution of (\ref{vtr22}). If $u\leq v$ on the parabolic boundary $B_R\times \{0\}\cup \partial B_R\times (0, T)$, then $u\leq v$ in $B_R\times (0, T)$.
\end{theorem}

The next theorem  from \cite[Theorem 5.3]{lara} is a parabolic version of   \cite[Proposition 2.9]{cx}.  It is the basic stability result which is needed in  compactness arguments.
\begin{theorem}\label{th4}(Stability)
	Let $\{\mathcal{F}_k\}_{k\geq 1}$ be a sequence of uniformly elliptic operators satisfying  (\ref{add1}).  Let $\{u^{(k)}\}_{k\geq 1}\subset C(B_R\times (-T, 0])$ be  viscosity solutions of \begin{equation*}\label{dd22}
		\mathcal{F}_k(D^2u^{(k)})\geq u^{(k)}_t\quad \text{in}\quad B_R\times (-T, 0].
	\end{equation*}
Assume that $\mathcal{F}_k\rightarrow \mathcal{F}$ uniformly in compact  subset  matrices $ \tilde{\mathcal{S} }\subset \mathcal{S},$
and $u^{(k)}\rightarrow u$ uniformly in compact subsets $E\subset B_R\times (-T, 0]$. Then $\mathcal{F}(D^2u)\geq u_t$ in the viscosity sense in $B_R\times (-T, 0]$.	
\end{theorem}

To summarize, we would like to mention that each of the above  results in this section will be
utilized  in proof of Theorem \ref{t1}. In particular,  the Harnack inequalities in Theorems \ref{th1} and \ref{thw1}  will be applied to prove the oscillation decay of the solutions to (\ref{a1}).  The regularity estimates in  Theorem \ref{th3} and comparision principle in Theorem \ref{th2}  are the key ingredients in our proof of the H$\ddot{\text{o}}$lder gradient estimates. Lastly the stability property in Theorem \ref{th4} is crucial to the proof of compactness and improvement of flatness.

\section{Fixed boundary problem}	
In this section,  we will  reduce the nolinear Stefan problem (\ref{a1}) into an equivalent nonlinear problem having  fixed boundary to study.  \vspace{0.2cm}

To do so, we first introduce a few notations useful to the subsequent clarification. Recall that
$B_r(x):=\{y\in \mathbb{R}^n: \,|x-y|<r\}$ are the balls in the $n$-dimensional Euclidean space. If $x$  is the origin, we will simply write $B_r$. Throughout the paper,  we usually write  $x=(x', x_n)\in  \mathbb{R}^{n-1}\times \mathbb{R}$, where $x_n=x\cdot e_n$ and  $e_n$ denotes the $n$-th  vector of the canonical basis of $\mathbb{R}^n$. Denote 
$$K_r:=(B_r(x)\cap \{x_n>0\})\times (-r,0],\quad F_r:=(B_r(x)\cap \{x_n=0\})\times (-r,0].$$
Define  the cubes $K_r^+$ and the Dirichlet boundary of $K_r$ respectively as
\begin{equation}\label{aer5}
K_r^+:= K_r\cup F_r,\quad \quad 	\partial_b K_r:=\partial K_r\cap (\{t=-1\}\cup \{x_n=1\}\cup_{i=1}^{n-1}\{|x_i|=1\}).
\end{equation}
Moreover,  denote a point in $\mathbb{R}^{n+1}$ as $X=(x,t)=(x', x_n, t)\in\mathbb{R}^{n-1}\times \mathbb{R}\times \mathbb{R}$.  We introduce the  distance between  $X_1(x,t)$ and $X_2(\tilde{x},\tilde{t})$
 \begin{equation}\label{a5}
 	d(X_1,X_2)=|x'-\tilde{x}'|+\min \left\{
 	|x_n|+|\tilde{x}_n|+|t-\tilde{t}|, \;|x_n-\tilde{x}_n|+|t-\tilde{t}|^{1/2}\right\}.
 \end{equation}
 From the above definition, it is clear to see that on the hyperplane $x_n=0$, the distance $d(X_1,X_2)$ is given by the standard Euclidean distance. If the points are far away from this hyperplane, we will adopt the parabolic distance. 
This distance in $x_n>0$ is naturally defined so that the rescaling property $|(\lambda x, \lambda^2 t)|=\lambda |(x,t)|$ holds.  In what follows, we
often use $c, c_i$ to denote small universal constants, and $C, C_i$ to denote large universal constants.

\subsection{Hodograph transform}
The Hodograph transform  is a helpful tool in free boundary problems. A well known application of this transform is in the higher regularity theory for $C^1$ free boundaries by Kinderlehrer-Nirenberg \cite{da}, and more recently for lower dimensional obstacle problems by Koch-Petrosyan-Shi \cite{he}.
\vspace{0.2cm}

Now we reformulate the problem (\ref{a1}) via the transformation of variable. More precisely, let  $u: \mathcal{Q}_{\lambda}\subset \mathbb{R}^{n+1}(n\geq 1)\rightarrow \mathbb{R}^+$  be the viscosity solution to (\ref{a1}). Suppose 
 \begin{equation}\label{a6}
 	y_{n+1}=u(y_1, \ldots, y_{n-1}, y_n, t)
 \end{equation}
and $\partial u/\partial y_n\not=0$. By  invert the relations, one has
 \begin{equation}\label{a7}
 		\begin{split}
 	y_{n}&=\tilde{u}(y_1, \ldots, y_{n-1}, y_{n+1}, t)\\[0.5ex]
 	&=:\tilde{u}(x_1, \ldots, x_{n-1}, x_{n}, t).
 		\end{split}
 	 \end{equation}
We notice that the graph  $\Gamma=\{(y,y_{n+1}, t)\,|\, y_{n}=\tilde{u}(y_1, \ldots, y_{n-1}, y_{n+1}, t)\}$ of $\tilde{u}$ is closed in $\mathbb{R}^{n+2}$, since it is  obtained as a rigid motion from the graph  $\Gamma=\{(y,y_{n+1}, t)\,|\, y_{n+1}=u(y_1, \ldots, y_{n}, t)\}$ of $u$.   And $\tilde{u}\subset  \mathcal{Q}_{\lambda}$ is a possibly muti-valued funtion  concerning $y_n$ direction.  \vspace{0.2cm}

Then the derivatives of $u$ may be expressed in terms of the derivatives of $\tilde{u}$ so that the equation (\ref{a1}) is transformed into one for $\tilde{u}$. From (\ref{a6}) and (\ref{a7}), it is not hard to deduce  
 \begin{equation}\label{a8}
 	u_t=-\frac{\tilde{u}_t}{\tilde{u}_{x_n}}\quad\text{and}\quad \nabla u=-\frac{1}{\tilde{u}_{x_n}}(\tilde{u}_{x_1},\ldots, \tilde{u}_{x_{n-1}}, -1).
\end{equation}
where $\nabla u$ denotes the gradient of $u$ with respect to the spacial variables only, $\tilde{u}_{t}:=\partial \tilde{u}/\partial t$ and  $\tilde{u}_{x_i}:=\partial \tilde{u}/\partial x_i$ for $i=1,\ldots, n$. From above we obtain expressions for the Hessian  of $u$,
$$D^2u=-\frac{1}{\tilde{u}_{x_n}} (A(\nabla \tilde{u} ))^TD^2 \tilde{u}\,A(\nabla \tilde{u}).$$
Here the matrix  $A\in \mathbb{R}^{n\times n}$ is given by

\renewcommand{\arraystretch}{1.8}
\[\left(\begin{array}{c|lr}
	I&\;\; 0   \\  [0.2ex]  \hline
-\frac{\tilde{u}_{x_1}}{\tilde{u}_{x_n}},\ldots, -\frac{\tilde{u}_{x_{n-1}}}{\tilde{u}_{x_n}}&\frac{1}{\tilde{u}_{x_n}} \\
	\end{array}\right)
\]
where $I$ is $(n-1)\times (n-1)$ identity matrix. The nonlinear Stefan problem (\ref{a1}) becomes
\be\label{a9}\le\{\begin{array}{lll}
	\tilde{u}_t=\widetilde{\mathcal {F}} (U)\quad\;\text{in}\quad K_r,\\[-0.5ex]
	\tilde{u}_t=\widetilde{\mathcal{G}}(\nabla \tilde{u})\quad \text{on}\quad F_r,
\end{array}\ri.\ee
where $\widetilde{\mathcal{G}}$ is some $C^2$ function,
\begin{equation}\label{aw10}
U=-\frac{1}{\tilde{u}_{x_n}} (A(\nabla \tilde{u} ))^TD^2 \tilde{u}\,A(\nabla \tilde{u})\quad \text{and}\quad \widetilde{\mathcal {F}}=-\tilde{u}_{x_n} \mathcal {F}.
 \end{equation}
Meanwhile,  for constant $M\geq 1$,  we assume
  \begin{equation}\label{a10}
 \nabla \tilde{u}\in B_{M}\cap \{\tilde{u}_{x_n}\geq M^{-1}\}=:\mathcal{O}_{M}.
 \end{equation}
By choosing $M$  large enough, we also suppose
 \begin{equation}\label{asw10}
 	\partial \widetilde{\mathcal{G}}/\partial (\nabla \tilde{u})_n\geq M^{-1}\quad \text{and} \quad \|\widetilde{\mathcal{G}}\|_{C^1}\leq M,
 	 \end{equation} 
  where $(\nabla \tilde{u})_n$ stands for the $n$-th variable of $\nabla \tilde{u}$. Moreover, $\tilde{u}$ solves  $(\ref{a9})$ in the viscosity sense as below. For completeness, we first give the notion of contact for the multi-valued function.

\begin{definition}\label{de4} 
We say that a single-valued function $\phi$ touches a multi-valued function $\tilde{u}$ above (resp. below)  at $(x_0, t_0)\in K_r\cup F_r$ in a parabolic cylinder $B_r(x_0)\times (t_0-r^2,t_0]$, if $\phi(x_0,t_0)\in \tilde{u}(x_0,t_0)$ and 
$$\tilde{u}(x,t)\leq \phi(x,t)\quad  (\text{resp.}\; u(x,t)\geq \phi(x,t))$$
for all possible values of $\tilde{u}$ at $(x,t)$ and for all $(x,t)\in B_r(x_0)\times (t_0-r^2,t_0]$.
\end{definition}

The definition of multi-valued viscosity solutions to (\ref{a9}) is formulated  as below. 
\begin{definition}\label{dek4} 
Assume that $\tilde{u}: \bar{K}_r\rightarrow \mathbb{R}$ is a multi-valued function with compact graph in $\mathbb{R}^{n+2}$. We say that $\tilde{u}$ is a viscosity subsolution (resp. supersolution)  to (\ref{a9})  if $\tilde{u}$ can not be touched by above  (resp. below)  at points in $K_r^+$ locally in parabolic cylinders by single-valued classical strict supersolutions (resp. subsolutions) $\phi$ of (\ref{a9}). The function $\tilde{u}$ is called a viscosity
solution if it is both a viscosity supersolution and a viscosity subsolution.
\end{definition}

One can see that $\tilde{u}$ is well defined in $B_{\bar{\lambda}}\times (-\bar{\lambda},0]$ with $\bar{\lambda}:=c\lambda$ for  $c$ small.  The $\bar{\epsilon}$-flatness assumption in Theorem \ref{t1} is written for $\tilde{u}$ (where we denote $\tilde{\lambda}$ by $\lambda$ for simplicity of notation)
	\begin{equation}\label{ass11}
	|\tilde{u}-l_{\tilde{a},\tilde{b}}(x,t)|\leq \epsilon \lambda\quad \text{in}\;\; K^+_{\lambda},
\end{equation}
where $l_{\tilde{a},\tilde{b}}=\tilde{a}_n(t) x_n+\tilde{b}(t)$ is defined as in (\ref{a3}) with 	
\begin{equation}\label{aa16}
	\tilde{b}'(t)=\widetilde{\mathcal{G}}(\tilde{a}_n(t)e_n).
		\end{equation}
We now set some conditions that will be used in the following proofs.  Assume 
	 \begin{equation}\label{af16}
	 0<\lambda\leq \min\{\lambda_0, \delta \epsilon\}
	 	\end{equation}
for small  and universal constants $\lambda_0, \delta>0$.  Assume in addition that  
	\begin{equation}\label{a16}
\tilde{a}(t)\in \mathcal{O}_M,\quad		|\tilde{a}^{\prime}_n(t)|\leq \delta \epsilon \lambda^{-2} \quad \text{with} \quad \epsilon\leq \epsilon_0,
	\end{equation}
 where the constant $\epsilon_0>0$ is universal. Owing  to (\ref{asw10}) and  (\ref{aa16})-(\ref{a16}), one has for $(x,t)\in (B_{\lambda}\cap \{x_n\geq 0\})\times  [t_0-\lambda^2, t_0+\lambda^2]\subset \overline{K}_{\lambda}$, 
 \begin{equation}\label{aff16}
 	\begin{split}
 	|l_{\tilde{a},\tilde{b}}(x,t)-\tilde{a}(t_0)\cdot x-\tilde{b}(t_0)| &\leq |\tilde{a}-\tilde{a}(t_0)|\cdot x+|\tilde{b}-\tilde{b}(t_0)|\\[1.5ex]
 	&\leq  C\delta \epsilon \lambda+C\lambda^2\leq C\delta \epsilon \lambda,
 	\end{split}
 \end{equation}
where $C>0$ is a universal constant. This together with (\ref{ass11}) leads to
\begin{equation}\label{a17}
	|\tilde{u}- \tilde{a}(t_0)\cdot x-\tilde{b}(t_0)|\leq C\epsilon \lambda,
\end{equation}
where $ a_n(t_0)\in \mathcal{O}_M$. Recalling (\ref{a7}), we see that  $\tilde{u}$ is in fact graphical in the $e_n$ direction. Together with (\ref{a17}), it implies that $\tilde{u}$ is single-valued in the region $x_n\geq C\epsilon \lambda$. Thus by choosing  $\epsilon_0$ small, for $r\in [\epsilon^{1/2}\lambda, \lambda]$, one can use the standard Harnack inequality (see e.g. \cite[Theorem 4.18]{wlh1}) to the solution $\tilde{u}$ of  the interior  equation  in  (\ref{a9}). \vspace{0.2cm}
 
 In order to improve the flatness of $\tilde{u}$, the general strategy we follow is to  linearize the equation near $l_{\tilde{a},\tilde{b}}$. For this purpose,  we define an error term for $(x,t)\in \overline{K}_{\lambda}$, 
\begin{equation}\label{a13}
	\vartheta\left(\frac{x}{\lambda}, \frac{t}{\lambda} \right):=\frac{1}{\epsilon \lambda} \left(\tilde{u}(x,t)-l_{\tilde{a},\tilde{b}}(x,t)\right).
\end{equation}
 In light of  (\ref{a9}) and (\ref{a13}), one can check  that $\vartheta$ solves
\be\label{a14}\le\{\begin{array}{lll}
&\lambda \vartheta_t=-\frac{\lambda \varTheta_n}{\epsilon} \mathcal {F} \left(-\frac{\epsilon}{\lambda \varTheta_n}  (A(\tilde{a}_n(\lambda t) e_n+\epsilon \nabla \vartheta ))^T D^2 \vartheta \,A(\tilde{a}_n(\lambda t)e_n+\epsilon \nabla \vartheta)\right)\\
&\quad \quad \quad -\frac{\lambda^2}{\epsilon} \tilde{a}^{\prime}_n(\lambda t) x_n-\frac{\lambda}{\epsilon}\tilde {b}^{\prime}(\lambda t)\quad \quad \quad \quad\; \,\text{in}\quad K_1,\\
& \; \vartheta_t=\frac{1}{\epsilon}\widetilde{\mathcal{G}}(\tilde{a}_n(\lambda t) e_n+\epsilon \nabla \vartheta)-	\frac{1}{\epsilon}\tilde{b}^{\prime}(\lambda t)\quad \quad  \text{on}\quad  F_1,
\end{array}\ri.\ee
where $M^{-1}\leq |A^T|, |A|\leq M$ and $\varTheta_n:= \tilde{a}_n(\lambda t)+\epsilon \vartheta_{x_n}=\tilde{u}_{x_n}\in \mathcal{O}_M$.

\subsection{Harnack inequality for $\vartheta$}
In this section we show that if a solution $\tilde{u}$ is sufficiently flat in a domain then the oscillation of $\vartheta$ of (\ref{a14}) decreases in a  smaller domain.  The notation $\mathscr{S}(\Lambda^{-1}, \Lambda)$ will stand for the class of all uniformly elliptic operators $\mathcal{F}$ with ellipticity constants $\Lambda$ and $\Lambda^{-1}$ with $\mathcal{F}(0)=0$. If $\mathcal{F}\in \mathscr{S}(\Lambda^{-1}, \Lambda)$, then 
\begin{equation}\label{ffd1}
	\mathcal{P}^{-}\Big(\mathcal{M}, \frac{1}{n\Lambda}, \Lambda\Big)\leq \mathcal{F}(\mathcal{M})\leq   \mathcal{P}^{+}\Big(\mathcal{M}, \frac{1}{n\Lambda}, \Lambda\Big)
\end{equation}
and the rescaled operator 
\begin{equation}\label{ff1}
	\mathcal{F}_\rho(\mathcal{M})=\frac{1}{\rho}	\mathcal{F}(\rho \mathcal{M})\in \mathscr{S}(\Lambda^{-1}, \Lambda), \quad \rho<0.
\end{equation}
In fact, for $\rho<0$ and any $\mathcal{N}\geq 0$, one has by (\ref{add1})  that
\begin{equation*}
	\begin{split}
\frac{1}{\rho}	\mathcal{F}(\rho (\mathcal{M}+\mathcal{N}))- \frac{1}{\rho}	\mathcal{F}(\rho\mathcal{M})&=-\frac{1}{\rho}\left( \mathcal{F}(\mathcal{A}-\rho \mathcal{N})-\mathcal{F}(\mathcal{A}) \right)\leq \Lambda \|\mathcal{N}\|		
	\end{split}
\end{equation*}
where $\mathcal{A}:=\rho (\mathcal{M}+\mathcal{N})$,  and vice versa.  Thus (\ref{ff1}) is true. 
Since $M^{-1}I\leq A:=A(\nabla \tilde{u})\leq M I$,  from (\ref{ffd1}), (\ref{ff1}) and Theorem \ref{thf1}, 
\begin{equation}\label{pt2}
		\begin{split}
	 \mathcal{P}^{-}_{\Lambda, M}(D^2 \vartheta)& \leq  \mathcal{P}^{-}\Big(A^TD^2 \vartheta A, \frac{1}{n\Lambda}, \Lambda\Big) \leq -\frac{\lambda \varTheta_n}{\epsilon}  \mathcal{F}\Big(-\frac{\epsilon}{\lambda \varTheta_n}  A^TD^2 \vartheta A\Big)\\
	 &\leq \mathcal{P}^{+}\Big(A^TD^2 \vartheta A, \frac{1}{n\Lambda}, \Lambda\Big) \leq \mathcal{P}^{+}_{\Lambda, M}(D^2 \vartheta),
	 	\end{split}
\end{equation}
where we write for simplicity
\begin{equation}\label{ptr2}
	\mathcal{P}^{\pm}_{\Lambda, M}(D^2 \vartheta):=\mathcal{P}^{\pm}\Big(D^2 \vartheta, \frac{1}{n\Lambda M^2}, \Lambda M^2\Big). 
\end{equation}

The first result, given below, is the interior Harnack inequality for $\vartheta$. In its proof, we will use the notations defined in (\ref{frq1}). 
\begin{lemma}\label{lem3}
Let $\vartheta$ be a viscosity solution to (\ref{a14}). Under the assumptions (\ref{aa16})-(\ref{a16}), if there is  $\mu\geq C\delta r^2$ such that
$$\varrho\leq \vartheta\leq\varrho+\mu \quad \text{in}\quad \mathcal{C}_r^{\lambda}\subset K_{1},$$
where the universal constants $\varrho, C>0$ ($C$ is large) and $r\geq \epsilon^{1/2}$, then  at least one of the following  holds in $\mathcal{C}_{r/2}^{\lambda}$,
$$either\;\;\quad  \varrho+c_0\mu\leq \vartheta\leq \varrho+\mu, \quad\; or \;\;\quad  \varrho\leq \vartheta\leq \varrho+(1-c_0)\mu,$$
where  $0<c_0<1$ is constant.  
\end{lemma}
\proof Consider the cylinder $\mathcal{C}_{r\lambda}(x_0, t_0)\subset K_{\lambda}$. For $(x,t)\in \mathcal{C}_{r\lambda}(x_0, t_0)$, we know by (\ref{a13}) that $\vartheta$ is defined in $\mathcal{C}_r^{\lambda}\subset K_{1}$. Let
\begin{equation}\label{prep1}
	\beta\left(x, t\right):=\vartheta\left(x, \lambda t\right)-\varrho.
\end{equation}
From  (\ref{af16}), (\ref{a16}) and (\ref{pt2}),  one finds that $\beta$ solves
\be\label{af9}\le\{\begin{array}{lll}
\beta_t\leq \mathcal{P}^{+}_{\Lambda, M}(D^2 \beta)+\hat{C} \delta\quad \text{in}\quad \mathcal{C}_{r}(x_0, t_0),\\
\beta_t\geq \mathcal{P}^{-}_{\Lambda, M}(D^2 \beta)-\hat{C}\delta \quad \text{in}\quad \mathcal{C}_{r}(x_0, t_0),
\end{array}\ri.\ee
where $\hat{C}>0$ is some constant. 
Since $\beta\geq 0$, thanks to Theorem \ref{th1}, 
\begin{equation}\label{pr1}
	\mathop{\sup \;\; \beta}_{\mathcal{C}^{-}_{r/2}(x_0, t_0)}\leq \tilde{C} \Big( \mathop{\inf\;\;  \beta }_{\mathcal{C}^{+}_{r/2}(x_0, t_0)}+\hat{C}\delta r^2\Big)
\end{equation}
for universal constant $\tilde{C}$.
Then we split into two cases to discuss. \vspace{0.2cm}\\
(i) Supoose that $\beta (x_0,t_0)\geq \mu/2$. Since $\mu\geq 3\tilde{C}\hat{C} \delta r^2$, we deduce  by  (\ref{pr1}) that
$$  \mathop{\inf\;\;  \beta }_{\mathcal{C}^{+}_{r/2}(x_0, t_0)}\geq \frac{1}{\tilde{C}}\beta (x_0,t_0) -\hat{C}\delta r^2 \geq \frac{\mu}{8\tilde{C}}.$$
Combining this and   (\ref{prep1}) yields
\begin{equation*}\label{dp2}
	\vartheta\geq  c_0\mu+\varrho \quad \text{in}\quad \mathcal{C}_{r/2}^{\lambda},
\end{equation*} 
where $c_0=\frac{1}{8\tilde{C}}$, as desired. \vspace{0.2cm}\\
(ii) Supoose that $\beta (x_0,t_0)\leq \mu/2$. It is clear that
$$\tilde{\beta}\left(x,  t\right):=\mu-\beta\left(x,  t\right)=\mu+\varrho-\vartheta\left(x,  \lambda t\right)\geq 0$$
 and $\tilde{\beta}$ satisfies (\ref{af9}). Moreover, $\tilde{\beta} (x_0,t_0)\geq \mu/2$. By the same argument as above, we obtain that
$\tilde{\beta}$ satisfies (\ref{pr1}). Therefore, 
$$ \vartheta\leq \varrho+(1-c_0)\mu \quad \text{in}\quad \mathcal{C}_{r/2}^{\lambda}.$$
This ends the proof of the lemma. $\hfill\Box$\\

Next we shall establish the Harnack inequality at the boundary.  To proceed, considering  (\ref{a5}), we define the parabolic cubes centered at $(x, t)$ with radius $r$, 
	\begin{equation}\label{a24}
	Q_r(x,t):=\le\{\begin{array}{lll}
		B_r(x)\times (t-r^2, t)\quad\quad  \quad\quad\quad \quad\;\, \text{if}\;r<|x_n|,\\[-0.5ex]
		\left(B_r(x)\cap \{x_n\geq 0\}\right)\times (t-r, t)\quad\quad  \text{if}\;|x_n|\leq r\leq1.
	\end{array}\ri.\end{equation}	
For $0<\lambda\leq 1$, one may perform a dilation of factor $\lambda^{-1}$ that maps  the domain $K_\lambda$ into $K_1$. So we apply the following rescaled distance 
	\begin{equation}\label{aqr24}
d_{\lambda}(P_1,P_2):=\lambda^{-1}d(\lambda P_1, \lambda P_2).
\end{equation}	
Then the rescaled  balls  are
	\begin{equation}\label{a25}
	Q^{\lambda}_{r}(x,t):=\le\{\begin{array}{lll}
		B_r(x)\times (t-\lambda r^2, t)\quad\quad  \quad\quad\quad \quad \text{if}\;r<|x_n|,\\[-0.5ex]
		\left(B_r(x)\cap \{x_n\geq 0\}\right)\times (t-r, t)\quad\quad  \text{if}\;|x_n|\leq r\leq \lambda^{-1}.
	\end{array}\ri.\end{equation}

Let  $osc_{Q_r} \vartheta:= \sup_{Q_r} \vartheta-\inf_{Q_r} \vartheta$.  We now come to  the result of  partial Harnack inequality for the solution $\vartheta$ to (\ref{a14}).
\begin{lemma}\label{lem8}
Under the assumptions of (\ref{a10})-(\ref{a16}), for any $(x,t)\in K_{1/2}$, it holds 
	\begin{equation*}\label{a26}
osc_{Q^{\lambda}_{r}(x,t)\cap K_{1}^+}\vartheta \leq C r^\alpha,	\quad \alpha\in (0,1),
\end{equation*}
where $r\geq C(\delta)\epsilon^{1/2} $ and the constant $\delta>0$ is small and universal. 
\end{lemma}
\proof  To clarify, we shall decompose the proof into three steps. In steps 1 and 2, we prove the oscillation decay of $\vartheta$, and then we do iteration in  step 3 to conclude the lemma. \vspace{0.2cm}

{\it Step 1.} 
From (\ref{ass11}), it is obvious to see  $|\vartheta|\leq 1$ in $\overline{K}_1$.  We first prove 
	\begin{equation}\label{a27}
osc_{K^+_{1/2}}\vartheta \leq 2(1-c)
\end{equation}	
for a universal constant $0<c<1$.  To achieve  this, we work on the function 
\begin{equation*}\label{dfg2}
	\tilde{\vartheta}(x,t)= \vartheta(x,t)+C\delta(2+t-x_n^2)+1,
\end{equation*}
where $C$ is large enough. Note that $\tilde{\vartheta}\geq 0$. Since 
\begin{equation}\label{dfgd2}
	osc_{E}\vartheta\leq osc_{E}\tilde{\vartheta}+osc_{E} \,C\delta  (x_n^2-t),\quad E\subset K_1^+,
\end{equation}
it is enough to prove  an oscillation decay of $\tilde{\vartheta}.$ Suppose $\epsilon\leq \epsilon_1(\delta)$ where $ \epsilon_1$ is a small constant depending only on $\delta$.  
By (\ref{a14}), (\ref{pt2}) and our assumptions, 
\begin{equation}\label{dfg1}
	\begin{split}
		 \mathcal{P}^{-}_{\Lambda, M}(D^2 \tilde{\vartheta}) &\leq\mathcal{P}^{-}\Big(A^TD^2 \tilde{\vartheta }A, \frac{1}{n\Lambda}, \Lambda\Big) +\mathcal{P}^{-}\Big(2C\delta IA^TA, \frac{1}{n\Lambda}, \Lambda\Big)\\
		 &\leq\mathcal{P}^{-}\Big(A^TD^2 \vartheta A, \frac{1}{n\Lambda}, \Lambda\Big)\leq -\frac{\lambda \varTheta_n}{\epsilon}  \mathcal{F}\Big(-\frac{\epsilon}{\lambda \varTheta_n}  A^TD^2 \vartheta A\Big)\leq \lambda \tilde{\vartheta}_t,
	\end{split}
\end{equation}
where the first inequality is due to the fact that $A^TA$ is positive definite. Hence $\tilde{\vartheta}$ solves
	\begin{equation}\label{a30}
\le\{\begin{array}{lll}
\mathcal{P}^{-}_{\Lambda, M}(D^2 \tilde{\vartheta}) 	\leq \lambda \tilde{\vartheta}_t \quad\quad\quad \quad\quad\;\; \text{in}\quad K_1, \\[-0.5ex]
	\frac{1}{M}\tilde{\vartheta}_{x_n}^+-M\left(\tilde{\vartheta}_{x_n}^-+|\nabla_{x'} \tilde{\vartheta}|\right)\leq \tilde{\vartheta}_t \quad \text{on}\quad F_1
	\end{array}\ri.\end{equation}	
in the viscosity sense, where $a^+:=\max\{a, 0\}$ and $a^-=(-a)^+$.   Let a sequence of times
$$t_j:=-1+\lambda j\in \big[-1, -\frac{1}{2}\big),$$
where $j$ are non-negative integers.   We denote  the last value $\tilde{j}$ such that $t_{\tilde{j}}<-1/2$ and $\tilde{j}+1=N_1+ N_2$ where $N_1\geq N_2\geq 1$. 
One claims that if there are $N_1$ times of $j$ such that 
	\begin{equation}\label{a38}
\tilde{\vartheta}\left(\frac{e_n}{2}, t_j+\frac{\lambda}{4}\right)\geq 1,
\end{equation}	
then
\begin{equation}\label{a39}
\tilde{\vartheta}(x,t)\geq \tilde{c}\quad \text{for}\quad (x,t)\in K^+_{1/2},
\end{equation}
where $\tilde{c}$ is small and universal constant. This together with (\ref{dfgd2}) implies the desired result (\ref{a27}). 
Conversely,  for $N_1$ times of  $j$ such that 
$$\tilde{\vartheta}\left(\frac{e_n}{2}, t_j+\frac{\lambda}{4}\right)<1,$$
we define $\hat{\vartheta}=2-\tilde{\vartheta}$, where $\hat{\vartheta}$ solves
\begin{equation*}\label{add30}
	\le\{\begin{array}{lll}
		\lambda \hat{\vartheta}_t\leq \mathcal{P}^{+}_{\Lambda, M}(D^2 \hat{\vartheta}) \quad\quad\quad\quad\quad\;\;  \text{in}\quad K_1, \\[-0.5ex]
		\hat{\vartheta}_t \leq M\left(\hat{\vartheta}_{x_n}^++|\nabla_{x'} \hat{\vartheta}|\right)-	\frac{1}{M}\hat{\vartheta}_{x_n}^- \quad \text{on}\quad F_1.
	\end{array}\ri.\end{equation*}	
Then one has $\hat{\vartheta} \left(\frac{e_n}{2}, t_j+\frac{\lambda}{4}\right)> 1$. Applying (\ref{a39}) to $\hat{\vartheta}$ leads to
\begin{equation*}\label{aa39}
	\hat{\vartheta}=2-\tilde{\vartheta}\geq \tilde{c}\quad \text{for}\quad (x,t)\in K^+_{1/2}.
\end{equation*}
And (\ref{a27}) is thereby  confirmed. \vspace{0.2cm}

{\it Step 2.}  We next prove the claim (\ref{a39}).  Denote $\beta(x,t):=T(t)\, h(x)$, where $h(x): \overline{K}_1\rightarrow \mathbb{R}$ satisfies
\begin{equation}\label{a40}
	\le\{\begin{array}{lll}
\mathcal{P}^{-}_{\Lambda, M}(D^2 h)=0  \quad \text{in}\quad B_{\frac{3}{4}}	\cap \{x_n>0\},	\\[-0.5ex]
h(x)=1 \quad\quad\quad \;\;\;\text{on} \quad B_{\frac{1}{2}}	\cap \{x_n=0\},\\[-0.5ex]
0\leq h(x)\leq 1 \quad\quad\,\text{on} \quad (B_{\frac{3}{4}}\backslash B_{\frac{1}{2}})	\cap \{x_n=0\}, \\[-0.5ex]
h(x)=0 \quad \quad\quad\quad \text{in} \quad (\overline{B}_1\backslash  B_{\frac{3}{4}}) \cap \{x_n\geq0\},
\end{array}\ri.\end{equation}	
and  for some $t_{j_0}\in (-1, 0]$,\vspace{0.1cm}
\begin{equation}\label{atr40}
T(t)=e^{-C_0(t-t_{j_0})}T(t_{j_0})\geq 0
\end{equation}	
with the universal constant $C_0>0$ chosen large. Note that $0\leq h\leq 1$ in $\overline{B}_1 \cap \{x_n\geq0\}$ and on $\{x_n=0\}\cap \{h=0\}$, we have $h_{x_n}>0$ and $|\nabla_{x'} h|=0$. Owing to this and (\ref{atr40}), it is not hard to check that $\beta$ solves
\begin{equation}\label{ax30}
	\le\{\begin{array}{lll}
		\mathcal{P}^{-}_{\Lambda, M}(D^2 \beta) 	\geq \lambda \beta_t \quad\quad\quad\quad\quad\;\; \text{in}\quad K_1, \\
		\frac{1}{M}\beta_{x_n}^+-M\left(\beta_{x_n}^-+|\nabla_{x'} \beta|\right)\geq \beta_t \quad \text{on}\quad F_1.
	\end{array}\ri.\end{equation}	
Therefore, if in $B_1\cap \{x_n\geq 0\}$,
\begin{equation}\label{axw30}
	\tilde{\vartheta}(x, t_{j_0})\geq \beta(x,t_{j_0})=T(t_{j_0})h(x),
\end{equation}	
by comparision  principle,  for $t\in [t_{j_0}, 0]$, one has 
\begin{equation}\label{a50}
\tilde{\vartheta}(x, t)\geq \beta(x,t)\quad \text{in}\quad B_1\cap \{x_n\geq 0\}.
\end{equation}	
Moreover under the hypothesis (\ref{axw30}), if $T(t_{j_0})\leq c_0$ and (\ref{a38})  holds for $t_{j_0}$, we assert that
\begin{equation}\label{a51}
	\tilde{\vartheta}(x, t_{j_0+1})\geq \beta(x,t_{j+1})=T(t_{j_0+1}) h(x) 
\end{equation}	
with $T(t_{j_0+1})\geq T(t_{j_0})+c_0\lambda$. Here $c_0>0$ is small and universal constant that will  be determined later. Otherwise assuming  (\ref{axw30}) holds,  if the values $j_0$ not satisfying (\ref{a38}) or $T(t_{j_0})>c_0$, one may use (\ref{a50}) to obtain $	\tilde{\vartheta}(x, t_{j_0+1})\geq T(t_{j_0+1}) h(x) $, where   $t_{j_0+1}$ given by (\ref{atr40}) belongs to $[t_{j_0}, 0]$. \vspace{0.3cm}

Now we verify (\ref{a51}). It is sufficient to prove (\ref{a51}) is ture for the first value $j^{\ast}\geq0$ such that (\ref{a38})  holds.   To see this, we denote $$t^{\ast}_i:=t_{j^{\star}}+i \frac{\lambda}{4}\quad \text{where}\quad i\in [0, 1, \ldots, 4].$$
Since $t_{j^{\ast}}\in [-1, -\frac{1}{2})$, one has $t_{j^{\ast}+1}=t^{\ast}_4\in [-1, 0)$. Notice that  we have the assumptions
\begin{equation}\label{ag38}
	\tilde{\vartheta}\left(\frac{e_n}{2},t^{\ast}_1\right)\geq 1, \,\quad T(t^{\ast}_{0})\leq c_0,\,\quad\tilde{\vartheta}(x, t^{\ast}_0)\geq T(t^{\ast}_{0}) h(x).
\end{equation}	
This together with Lemma \ref{lem3} (provided that $\epsilon^{1/2}\leq \hat{c}$ with $\hat{c}$ small and universal) implies
\begin{equation}\label{a52}
\tilde{\vartheta}(x,t)\geq c_1 \quad \text{in}\quad B_{\frac{3}{4}-\kappa}	\cap \{x_n\geq \kappa\}\times [t^{\ast}_2,t^{\ast}_4]\subset K_1,
\end{equation}	
where  the constants $0<c_1<1$ and $0<\kappa<1/4$ are universal. For simplicity, denote the spatial domain
$$B^{\ast}_{\kappa}:=\left(B_{\frac{3}{4}}	\cap \{x_n>0\}\right)\backslash \left( B_{\frac{3}{4}-\kappa}	\cap \{x_n\geq \kappa\}\right).$$
Now we construct barriers on $B^{\ast}_{\kappa} \times [t^{\ast}_2,t^{\ast}_4]$ to compare with $\tilde{\vartheta}$. Define 
\begin{equation*}\label{ass55}
		b(x,t)=	
(T(t^{\ast}_{3})+c_2(t-t^{\ast}_3)) h(x)+\frac{c_1}{2} (e^{\zeta(x)-1}-2) \quad \text{for}\quad t\in [t^{\ast}_2,t^{\ast}_4], 
\end{equation*}	
where the constant $0<c_2\leq \frac{c_1}{2}<1$, $h(x)$ is defined  in  $(\ref{a40})$ and  $0\leq \zeta(x)\leq 1: \overline{B^{\ast}_{\kappa}}\rightarrow \mathbb{R}$ satisfies 
\begin{equation}\label{a55}
	\le\{\begin{array}{lll}
		\mathcal{P}^{-}_{\Lambda, M}(D^2 \zeta(x))\geq \frac{2e}{c_1}  \quad \text{in}\quad B^{\ast}_{\kappa},	\\[-0.5ex]
		\zeta=0,\; |\nabla \zeta|\geq 1 \quad \quad\;\,\text{on}\quad \partial (B_{\frac{3}{4}}	\cap \{x_n\geq 0\}),\\[-0.5ex]
		\zeta \leq 1  \quad \quad\quad\quad\quad\quad\;\,\text{on}\quad  \partial ( B_{\frac{3}{4}-\kappa}	\cap \{x_n\geq \kappa\}).
	\end{array}\ri.\end{equation}	
Then on $B^{\ast}_{\kappa} \times [t^{\ast}_2,t^{\ast}_4],$ there holds
\begin{equation}\label{a56}
	\lambda b_t\leq \lambda c_2h(x)\leq c_2 \leq \frac{c_1}{2} e^{\zeta(x)-1} \mathcal{P}^{-}_{\Lambda, M}(D^2 \zeta(x))\leq  \mathcal{P}^{-}_{\Lambda, M}(D^2 b).
\end{equation}	
Now we compare $\tilde{\vartheta}$ and $b$ on the boundary. Notice  from  (\ref{a50}) and  (\ref{ag38}) that $\tilde{\vartheta}(x, t)\geq T(t^{\ast}_3) h(x)$.  On $\partial (B_{\frac{3}{4}}	\cap \{x_n\geq 0\})$ or at $t=t^{\ast}_2$, one has
\begin{equation}\label{ap56}
	b(x,t)\leq \left(T(t^{\ast}_{3})+\frac{c_2\lambda}{4}\right) h(x)-\frac{c_1}{2}\leq T(t^{\ast}_3) h(x) \leq \tilde{\vartheta}(x, t). 
\end{equation}	
While on $\partial ( B_{\frac{3}{4}-\kappa}	\cap \{x_n\geq \kappa\})$,  due to (\ref{ag38}) and (\ref{a52}),
\begin{equation}\label{ad57}
	b \leq T(t^{\ast}_{3}) h(x) +\frac{c_2}{2}+\frac{c_1}{2}\leq c_0+\frac{3c_1}{4}\leq c_1\leq \tilde{\vartheta}
\end{equation}
by choosing $c_0\leq c_1/4$.  
In view of (\ref{a56})-(\ref{ad57}), utilizing Theorem \ref{th2}, we obtain $\tilde{\vartheta}\geq b$ on $B^{\ast}_{\kappa} \times (t^{\ast}_3,t^{\ast}_4]$. It then follows from (\ref{a52}) that in $K_{\frac{3}{4}}^+$,
\begin{equation*}\label{a59}
		\begin{split}
\tilde{\vartheta}(x, t^{\ast}_4)=\tilde{\vartheta}(x, t_{j^{\star}}+\lambda)&\geq \left(T(t^{\ast}_{3})+\frac{c_2 \lambda}{4}\right) h(x)+\frac{c_1}{2} e^{\zeta(x)-1}\\
&\geq \left(T(t^{\ast}_{0})e^{-C_0\frac{3\lambda}{4}} +\frac{c_2 \lambda}{4}\right) h(x)\\
&\geq \left(T(t^{\ast}_{0})-T(t^{\ast}_{0})\left(1-e^{-C_0\frac{3\lambda}{4}}\right) +\frac{c_2 \lambda}{4}\right) h(x)\\
&\geq \left(T(t^{\ast}_{0})+c_0\lambda\right) h(x),
	\end{split}
\end{equation*}
where $c_0>0$ is small. Thus (\ref{a51}) is proven. \vspace{0.1cm}

Finally, since $\tilde{\vartheta}\geq 0$, we can set the initial data $T(t_{0})=0$.  Employ the results (\ref{a50}) and  (\ref{a51}) repeately to the sequence of times $t_j$.  Together with (\ref{atr40}), we infer 
$T(t_{\tilde{j}})\geq \tilde{c}>0$ and thus  $\tilde{\vartheta}(x, t_{\tilde{j}})\geq T(t_{\tilde{j}})h(x)\geq \tilde{c}h(x),$ where  $t_{\tilde{j}}<-1/2$,  and  the constant  $\tilde{c}$ depends only on $c_0$ and $C_0$. Observe by (\ref{a40}) that $h(x)\geq c>0$ in $B_{1/2}\cap \{x_n\geq 0\}$. Consequently, we conclude from (\ref{a50}) that
$$\tilde{\vartheta}(x,t)\geq T(t)\, h(x) \geq \tilde{c}\quad \text{in} \quad K^+_{1/2}.$$
Therefore the claim (\ref{a39}) is proved and (\ref{a27}) follows immediately. \vspace{0.2cm}

{\it Step 3.} Utilizing  (\ref{a27}), we shall end the proof of the lemma by an  iteration argument. Precisely, from (\ref{ass11}), (\ref{a13}) and (\ref{a27}), 
	\begin{equation}\label{a61}
	|\tilde{u}-l_{\tilde{a},\tilde{b}}(x,t)|\leq 2(1-c)\epsilon \lambda\quad \text{in}\;\; K^+_{\lambda/2}.
\end{equation}
By (\ref{a52}), (\ref{a61}) and the hypothesis $\epsilon\leq \epsilon_1(\delta)$, one infers that for integers $k\geq 1$,  
$$2^{k}\epsilon^{1/2}\leq \hat{c}\quad \text{and}\quad 2^k(1-c)^k\epsilon \leq \epsilon_1(\delta)\quad \text{in}\quad  K^+_{\lambda/2^k}.$$
This implies that we can iterate finite $k$ times provided that
	\begin{equation}\label{a62}
		2^{k}\epsilon^{1/2}\leq  \epsilon_2(\delta),
	\end{equation}
where the constant $\epsilon_2$ depends only on $\delta$. Meanwhile, thanks to  (\ref{a27}), for any $r\leq 1$, 
	\begin{equation*}\label{as70}
		osc_{K^+_{1/2}}\vartheta(rx, rs)\leq (1-c)osc_{K^+_{1}}\vartheta(rx, rs).
	\end{equation*}
Recall (\ref{a25}), and thus the above inequality implies
	\begin{equation}\label{ase70}
		osc_{Q^{\lambda}_{r/2}(0,0)}\vartheta(x, s)\leq (1-c)osc_{Q^{\lambda}_{r}(0,0)}\vartheta(x, s),\;\quad \forall r\leq 1.
\end{equation}
We claim that for every $(\hat{x},\hat{s})\in K_{1/2}$, there holds 
	\begin{equation}\label{a70}
osc_{Q^{\lambda}_{r/8}(\hat{x},\hat{s})}\vartheta(x, s)\leq (1-c)osc_{Q^{\lambda}_{r}(\hat{x},\hat{s})}\vartheta(x, s)\;\quad \text{for}\;\, C(\delta)\epsilon^{1/2} \leq  r\leq 1/4.
\end{equation}
To prove the assertion, we distinguish into two cases. If  $C(\delta)\epsilon^{1/2} \leq r< \hat{x}_n$, the claim  follows from Lemma \ref{lem3} and (\ref{a62}) directly.  If  $\hat{x}_n\leq r\leq 1/4$, for the boundary points $((\hat{x}',0), \hat{s})\in F_{1/2}$, one has by (\ref{ase70}) that
	\begin{equation*}\label{ad70}
	osc_{Q^{\lambda}_{r/4}((\hat{x}',0), \hat{s})}\vartheta \leq (1-c)osc_{Q^{\lambda}_{r/2}((\hat{x}',0). \hat{s})}\vartheta.
\end{equation*}
 Thus (\ref{a70}) is proved. Then let $(x,s)\in Q^{\lambda}_{r}(\hat{x},\hat{s})\cap K_1^+$ and  let $m>0$ be integers such that $$(x,s)\in Q^{\lambda}_{2^{-m}}(\hat{x},\hat{s})\backslash Q^{\lambda}_{2^{-m-1}}(\hat{x},\hat{s}).$$
 As a consequence of  (\ref{a70}),
	\begin{equation}\label{a63}
		|\vartheta(x,s)-\vartheta(\hat{x},\hat{s})|\leq osc_{Q^{\lambda}_{2^{-m}}(\hat{x},\hat{s})} \vartheta \leq C_02(1-c)^m=C_02^{1-\alpha m},
	\end{equation}
where $\alpha=-\log_2(1-c)$. Since $2^{-m-1}\leq |(x,s)-(\hat{x},\hat{s})|$, we derive from  (\ref{a63}) that
	\begin{equation*}\label{as63}
		osc_{Q^{\lambda}_r(x_0, s_0)\cap K_1^+}\vartheta=|\vartheta(x,s)-\vartheta(x_0,s_0)|\leq Cr^{\alpha}, 
	\end{equation*}
where $0<\alpha<1$ and  $r\geq C(\delta)\epsilon^{1/2}$. This  completes the proof of the lemma.  $\hfill\Box$

\subsection{The linearized problem}
This section is devoted to linearizing the nonlinear problem (\ref{a9}). 
Since $\lambda/\epsilon\leq \delta $, by  (\ref{aa16}), (\ref{a16}) and (\ref{a14}), letting $\delta\rightarrow0$ and $\epsilon\rightarrow0$,  we arrive at the linear equation
\be\label{a18}\le\{\begin{array}{lll}
	\lambda \eta_t=\mathcal {L}(A^{T}(\tilde{a}_n(\lambda t) e_n)D^2 \eta A(\tilde{a}_n(\lambda t)e_n)\quad \text{in}\quad K_1,\\[-0.5ex]
	\eta_t=\nabla\widetilde{\mathcal{G}}(\tilde{a}_n(\lambda t) e_n) \nabla \eta\quad \quad\quad\quad\quad\quad\quad\;\,\text{on}\quad F_1.
\end{array}\ri.\ee
Here $\mathcal {L}$ is a linear operator, precisely, 
\begin{equation*}\label{dfbc13}
	\begin{split}
		\mathcal {L}(A^T(\lambda t) D^2\eta A(\lambda t))=tr(PA^T(\lambda t) D^2\eta A(\lambda t))
		\quad \text{with}\quad \Lambda^{-1}I \leq P\leq \Lambda I,
			\end{split}
	\end{equation*} 
where  $A(\lambda t):=A(\tilde{a}_n(\lambda t)e_n)$, and  $\widetilde{\mathcal{G}}$ satisfies  (\ref{asw10}). Due to (\ref{a16}), 
	\begin{equation}\label{aw63}
		\;\quad |(A(\lambda t)PA^T(\lambda t))^{\prime}|\leq C\lambda^{-1},\quad  |(\nabla\widetilde{\mathcal{G}})^{\prime}|\leq \lambda^{-1}.
\end{equation}

Firstly, we state the central property of viscosity solutions to (\ref{a18}), i.e. comparison principle,  that will be used later. For  its proof we refer to \cite[Lemma 9.1]{r1} that can be adapted to our setting without difficulty. 
\begin{lemma}\label{lm1}
Let $\eta^{sub}$ be a viscosity subsolution and $\eta^{sup}$ be  a viscosity supersolution of (\ref{a18}) in $K_1^+$. If $\eta^{sub}\leq \eta^{sup}$ on $\partial_b K_1$, then $\eta^{sub}\leq \eta^{sup}$ in $K_1$,
where  $\partial_b K_1$ is defined as in (\ref{aer5}).
\end{lemma}

Let $E$ be an open subset of $\overline{K}_1$. For $\alpha\in (0,1)$, we next define the parabolic H$\ddot{\text{o}}$lder  norms and semi-norms 
\begin{equation}\label{re2}
	[\eta]_{C^{0,\alpha}_{x,\,t}(E)}:=\sup_{(x,t), (y,s)\in E, \atop (x,t)\not=(y,s)} \frac{|\eta(x,t)-\eta(y,s)|}{d_{\lambda} ((x,t),(y,s))^\alpha}, \quad\quad [\eta]_{C^{\alpha}_{t}(E)}:=\sup_{(x,t), (x,s)\in E, \atop t\not=s} \frac{|\eta(x,t)-\eta(x,s)|}{|t-s|^{\alpha}},
\end{equation}
where $d_\lambda$ is defined as in (\ref{aqr24}). Furthermore, we set
$$[\eta]_{C^{1,\alpha}_{x,\,t}(E)}:=[\nabla_x \eta]_{C^{0,\alpha}_{x,\,t}(E)}+[\eta]_{C^{\frac{1+\alpha}{2}}_{t}(E)}.$$
 For indexes $k\leq 1$, we say that $\eta\in C^{k,\alpha}_{x,\,t}(E)$ when $[\eta]_{C^{k,\alpha}_{x,\,t}(E)}<\infty$, where
$$\|\eta\|_{C^{k,\alpha}_{x,\,t}(E)}:=\sum_{i\leq k}\|D^i \eta\|_{L^\infty(E)}+[\eta]_{C^{k,\alpha}_{x,t}(E)}.$$
On the other hand, denote two classes of functions
$$\overline {\mathcal{S}}_{\Lambda, M}:=\Big\{ \eta\in C(K_1^+):\mathcal{P}^{-}_{\Lambda, M}(D^2 \eta) 	\leq \lambda \eta_t  \; \text{in}\; K_1,\; \frac{1}{M}\tilde{\eta}_{x_n}^+-M\left(\eta_{x_n}^-+|\nabla_{x'} \eta|\right) \leq  \eta_t\; \text{on}\; F_1  \Big\},$$ 
$$\underline {\mathcal{S}}_{\Lambda, M}:=\Big\{ \eta\in C(K_1^+):\mathcal{P}^{+}_{\Lambda, M}(D^2 \eta) 	\geq \lambda \eta_t  \; \text{in}\; K_1,\;  M\left(\eta_{x_n}^++|\nabla_{x'} \eta|\right)-	\frac{1}{M}\eta_{x_n}^-\geq \eta_t\; \text{on}\; F_1  \Big\}.$$ 

Then we investigate the  H$\ddot{\text{o}}$lder regularity for $\eta$, which is presented as below. \vspace{0.2cm} 
\begin{lemma}\label{lem9}
	Let $\eta\geq 0$ be a viscosity solution to (\ref{a18}). One has\vspace{0.1cm}  \\
	(i) there exists a universal constant $0<\mu<1$ such that
	\begin{equation*}\label{ad28}
		\text{osc}_{K^+_{1/2}}\eta \leq (1-\mu) osc_{K^+_1}\eta.
	\end{equation*}
	(ii) there  exists a H$\ddot{o}$lder exponent $0<\alpha<1$ and universal constant $C$ such that
	\begin{equation*}\label{af30}
		\|\eta\|_{C^{0,\alpha}_{x,\,t}(K^+_{1/2})}\leq C\|\eta\|_{L^\infty(K^+_1)}.
	\end{equation*}
	(iii) if in addition $\eta\in C(\partial_b K_1)$ and $\eta|_{\partial_b K_1}=\psi$ with  $\psi\in C^{0,\alpha}_{x,\,t}(\partial_b K_1)$, then 
	\begin{equation}\label{a31}
		\eta\in C^{0,\alpha}_{x,\,t}(\overline{K}_1)\quad \text{and}\quad \|\eta\|_{C^{0,\alpha}_{x,\,t}(\overline{K}_1)}\leq C \|\psi\|_{C^{0,\alpha}_{x,\,t} (\partial_b K_1)}.
	\end{equation}
\end{lemma}
\proof For the statements  (i) and (ii), note by (\ref{a18}) that $\eta\in \mathcal{S}_{\Lambda, M}:=\overline {\mathcal{S}}_{\Lambda, M}\cap \underline {\mathcal{S}}_{\Lambda, M}.$ In this setting, the subsequent proof is essentially identical to that of Lemma \ref{lem8}, so we omit the details. We now  prove (iii). The proof is divided into two steps. \vspace{0.2cm}
 
{\it{Step 1.}} From (i), one can in fact deduce that for  $(x,t)\in F_1$, if $ Q_{r}^{\lambda}(x,t)\cap \partial_b K_1=\emptyset$ where $r\leq 1/\lambda$, then 
\begin{equation*}\label{ds3}
			osc_{Q_{r/2}^{\lambda}(x,\,t)}\eta \leq (1-\mu) osc_{Q_{r}^{\lambda}(x,\,t)}\eta.
\end{equation*}
If $ Q_{r}^{\lambda}(x,t)\cap \partial_b K_1\not=\emptyset$, we claim that
\begin{equation}\label{d3}
	osc_{\overline{K}_1\cap Q_{r/2}^{\lambda}(x,\,t)}\eta \leq (1-\mu)  osc_{\overline{K}_1\cap Q_{r}^{\lambda}(x,\,t)}\eta.
\end{equation}
To see this, we define $V=\cup_{i=1}^{n}(-\infty, v_i)\times (-t, +\infty)$, where $v_i, t\in [0,1]$ and $\min_{1\leq i\leq n} v_i\leq 3/4$. Assume that $\eta\geq 0$ solves
\be\label{d18}\le\{\begin{array}{lll}
\eta \in \overline {\mathcal{S}}_{\Lambda, M} \quad \text{in} \quad K_1^+\cap V, \\[-0.5ex]
\eta\geq \frac{1}{6}  \quad \quad\;\,\text{in} \quad K_1^+\cap \partial V.
\end{array}\ri.\ee
For the converse, one may discuss for $\tilde{\eta}=1-\eta$ by a similar argument of the problem (\ref{d18}).   
Let $\check{\eta}=\min \{\eta, 1/6\}$  extended by $1/6$ in $K_1\backslash V$ and obviously $\check{\eta}\in  \overline {\mathcal{S}}_{\Lambda, M}$ in $K_1^+$. If $t<1$, by (\ref{a50}), there is a universal constant $c$ such that
\begin{equation}\label{da3}
	\eta\geq \check{\eta}\geq  c\quad \text{in}\quad K_{1/2}\cap V.
\end{equation}
This leads to (\ref{d3}), as claimed. If $t= 1$, since $\min_{1\leq i\leq n} v_i\leq 3/4$, for any $t\in [-1, 0]$,  there is a universal constant $c_0>0$ such that 
\begin{equation}\label{da4}
\left|\left\{ (x,t)\in K_1: \check{\eta}\geq 1/6\right\}\right|\geq c_0 |K_1|.
\end{equation}
Now choose $\rho>0$ small enough such that $1/\rho^2\in \mathbb{N}^+$ (depends on $c_0$ and $n$). Denote
$$E:=\{(x,t)\in E_{10\rho} \times (-1, -10\rho^2]\}\subset K_1, \quad E_{10\rho}:= dist(x, \partial (B_1\cap \{x_n>0\}))> 10\rho,$$
then we deduce
\begin{equation}\label{da7}
	\begin{split}
\left|\left\{ (x,t)\in E: \check{\eta}\geq 1/6\right\}\right|&\geq \left|\left\{ (x,t)\in K_1: \check{\eta}\geq 1/6\right\}\right|-|K_1\backslash E|\\
&\geq c_0 |K_1|- C(n)\rho\geq \frac{c_0}{2} |K_1|,
	\end{split}
\end{equation}
where $C(n)>0$ is a constant depending on $n$.  Define a sequence of  cylinders 
$$E^k=(B_{\rho}(x^{k})\cap \{x^{k}_n>\rho\})\times (-1+ (k-1)\rho^2, -1+ k\rho^2],$$
where $k=1,2, 3,\cdots$, and $x^k\in \overline{E}_{10\rho}.$ By moving the slices one by one, we find that there are finite $N\geq 1/\rho^2-10$ cylinders satisfying $E\subset \cup_{k=1}^{N}E^k$.  So there
exists at least one cylinder $E^{k_0}$ with $x^{k_0}\in E_{11\rho}$ such that
\begin{equation}\label{das8}
	\left|\left\{ (x,t)\in E^{k_0}: \check{\eta}\geq 1/6\right\}\right|\geq  \frac{c_0}{2N} |K_1|.
\end{equation}
Suppose not, 
\begin{equation*}\label{da8}
	\begin{split}
		\left|\left\{ (x,t)\in E: \check{\eta}\geq 1/6\right\}\right|&\leq \left|\cup_{k=1}^N (x,t)\in E^k: \check{\eta}\geq 1/6 \right|\\
		&\leq \sum_{k=1}^N  \left|(x,t)\in E^k: \check{\eta}\geq 1/6 \right|<\frac{c_0}{2} |K_1|,
	\end{split}
\end{equation*}
contradicting to (\ref{da7}). It then follows from Theorem \ref{thw1}  and (\ref{das8}) that 
\begin{equation}\label{da5}
	\check{\eta}\geq c_1 \quad \text{in} \quad (B_{\rho}(x^{k_0})\cap \{x^{k_0}_n>0\})\times (-1+ \rho^2, -1+ 2\rho^2] \subset K_1,
\end{equation}
where the constant $c_1>0$ depending on $c_0$, $n, \Lambda$ and $M$.  Thus considering $t\in [t_0, t_0+\lambda/4]$ in (\ref{da4}), after proper rescaling, we can get the result as in (\ref{a52}), and the remaining proof is argued similarly as in  Lemma \ref{lem8}. 
Thus (\ref{d3}) is ture and  in conclusion 
\begin{equation}\label{d4}
		\|\eta\|_{C^{0,\alpha}_{x,t}(\bar{K}_{1}\cap Q_{r/2}^{\lambda}(x,\,t))}\leq C\|\eta\|_{L^\infty(\bar{K}_{1}\cap Q_{r}^{\lambda}(x,\,t))}\quad \forall r\leq 1/\lambda,
\end{equation}
where $(x,\, t)\in K^+_{1}$. \vspace{0.2cm}

{\it{Step 2.}} Under the hypothesises of (iii),  one can check that for any $(\hat{x},\hat{t})\in  \partial_b K_1$, 
\begin{equation}\label{d5}
\sup_{(y,s)\in \overline{K}_1, (y,s)\not=(\hat{x},\hat{t})} \frac{|\eta(y,s)-\eta(\hat{x},\hat{t})|}{d_{\lambda} ((y,s),(\hat{x},\hat{t}))^{\alpha}}\leq C \sup_{(\hat{y},\hat{s})\in \partial_b K_1}\frac{|\psi(\hat{y},\hat{s})-\psi(\hat{x},\hat{t}) |}{d_{\lambda} ((\hat{y},\hat{s}),(\hat{x},\hat{t}))^\alpha},
\end{equation}
where $C$ is universal constant. And by the maximum principle, 
\begin{equation}\label{d8}
	\inf_{\partial_b K_1}\psi\leq \eta\leq\sup_{\partial_b K_1}\psi \quad \text{in}\quad  K_1^+.
\end{equation}
Consider  two points  $(x,t)=(x', x_n, t), (y,s)=(y', y_n, s)\in K_1$.   Let 
$$d_{x,t}=\text{dist}((x,t), \partial_b K_1),\quad \quad d_{y,s}=\text{dist}((y,s), \partial_b K_1).$$ Assume $d_{x,t}\leq d_{y,s}$ and  $x_n\leq y_n$. Take $(\hat{x}, \hat{t}), (\hat{y}, \hat{s})\in \partial_b K_1$ such that
$$d_\lambda((x,t),(\hat{x}, \hat{t}))=d_{x,t},\quad \quad d_\lambda((y,s),(\hat{y}, \hat{s}))=d_{y,s}.$$
{\it{Case 1}}. Assume $d_\lambda((x,t),(y, s))\leq d_{y,s}/2$.  Then $(x,t)\in \bar{Q}^{\lambda}_{d_{y,s}/2}(y,s), $ where $Q^\lambda_r$ is defined as in (\ref{a25}). 
 Clearly $\bar{Q}^{\lambda}_{d_{y,s}/2}(y,s)\cap \partial_b K_1=\emptyset.$  In this case, there are two possibilities.  \vspace{0.2cm}
 
\noindent (i) If $ d_{y,s}/2>y_n$,  it then follows from  (\ref{d5}) and the interior estimates that
\begin{equation}\label{dff9}
	d_{y,s}^{\alpha}[\eta]_{C^{0,\alpha}_{x,\,t}(K_1)}\leq C \|\eta-\eta(\hat{y}, \hat{s})\|_{L^\infty(Q^{\lambda}_{d_{y,s}}(\hat{y}, \hat{s}))}\leq C d_{y,s}^{\alpha} \|\psi\|_{C^{0,\alpha}_{x,\,t} (\partial_b K_1)}.
\end{equation}

\noindent (ii) If $ d_{y,s}/2\leq y_n$,  from (\ref{d4}) and (\ref{d8}), we infer 
\begin{equation}\label{d9}
	\begin{split}
&\quad 	|\eta(x,t)-\eta(y,s)|\\
	&\leq |\eta(x,t)-\eta(x',0,t)|+|\eta(x',0,t)-\eta(y',0,s)|+|\eta(y',0,s)-\eta(y,s)|\\[0.5ex]
	&\leq C \left(d_{\lambda} ((x,t),(x',0,t))^\alpha + d_{\lambda} ((x',0,t),(y',0,s))^\alpha+d_{\lambda} ((y',0,s),(y,s))^\alpha \right) \|\eta\|_{L^\infty(K_1^+)}\\[0.5ex]
	&\leq C(x_n^\alpha+d_{\lambda} ((x',0,t),(y',0,s))^\alpha+y_n^\alpha) \|\eta\|_{L^\infty(K_1^+)}\\[0.5ex]
	&\leq C d_{\lambda} ((x,t),(y,s))^{\alpha} \|\psi\|_{C^{0,\alpha}_{x,\,t} (\partial_b K_1)}.
	\end{split}
\end{equation}
Thus (\ref{a31}) follows by (\ref{dff9}) and (\ref{d9})  immediately.\vspace{0.2cm}

\noindent {\it{Case 2}}. Assume $d_\lambda((x,t),(y, s))\geq d_{y,s}/2$. If $ d_{y,s}/2>y_n$, note that 
$$d_{\lambda} ((\hat{x},\hat{t}), (\hat{y},\hat{s}))\leq d_{x,t}+d_\lambda((x,t),(y, s))+ d_{y,s}\leq 5 d_\lambda((x,t),(y, s)). $$
Due to this and (\ref{d5}), there holds
\begin{equation}\label{d10}
	\begin{split}
		|\eta(x,t)-\eta(y,s)|&\leq |\eta(x,t)-\eta(\hat{x}, \hat{t})|+|\eta(\hat{x}, \hat{t})-\eta(\hat{y}, \hat{s})|+|\eta(\hat{y}, \hat{s})-\eta(y,s)|\\[0.5ex]
		&\leq C \left( d_{x,t} ^{\alpha} +d_{\lambda} ((\hat{x},\hat{t}), (\hat{y},\hat{s}))^{\alpha}  + d_{y,s}^{\alpha} \right) \|\psi\|_{C^{0,\alpha}_{x,\,t} (\partial_b K_1)}\\[0.5ex]
		&\leq Cd_\lambda((x,t),(y, s))^{\alpha}   \|\psi\|_{C^{0,\alpha}_{x,\,t} (\partial_b K_1)}. 
	\end{split}
\end{equation}
If $ d_{y,s}/2\leq y_n$, similarly as above, by  (\ref{d5}), (\ref{d9}) and (\ref{d10}), we arrive at (\ref{a31}) eventually. The proof of the lemma is finished. $\hfill\Box$

\subsection{$C^{1,\alpha}$ estimates and existence of  $\eta$.}
Relying on the $C^{0, \alpha}$ estimate for $\eta$ established above, we are now in a position to prove its $C^{1,\alpha}$ regularity, namely,  
\begin{lemma}\label{lem11}
	Let $\eta$ be a viscosity solution to (\ref{a18}). Assume $\|\eta\|_{L^\infty}\leq 1$. Then for universal constant $C>0$, one has\\
(i) 
\begin{equation*}\label{a33}
|\nabla \eta|\leq C\quad \text{and} \quad |D^2 \eta|\leq C\quad \forall (x,t)\in K_{1/2},
\end{equation*}	
(ii) for any $\rho\leq 1/2$, there are  $\bar{a}(t)$ and $\bar{b}(t)$ satisfying 
\begin{equation}\label{awq34}
	|\bar{a}|\leq C, \quad |\bar{a}^{\prime}_n(t)|\leq C \rho^{\alpha-1}\lambda^{-1}\quad \text{and}\quad \bar{b}^{\prime}(t)=\nabla\widetilde{\mathcal{G}}(\tilde{a}_n(\lambda t) e_n) 
\end{equation}
such that in $K_{\rho}$,
\begin{equation*}\label{a34}
	|\eta-\bar{a}(t)x-\bar{b}(t)|\leq C\rho^{1+\alpha},\quad  0<\alpha<1.
\end{equation*}	
\end{lemma}
\proof The proof is split into three steps. The initial step is to  complete (i) by establishing H$\ddot{\text{o}}$lder estimates for solutions to 1-D linear problem. Based on this, we then conclude (ii) via modifying the linear approximation.  \vspace{0.2cm}

\emph{Step 1}. Note that $osc_{Q_r^\lambda} \eta=sup_{Q_r^\lambda}\eta-inf_{Q_r^\lambda}\eta$. By adding or subtracting an appropriate
constant,  we deduce from Theorem \ref{th3} and Lemma \ref{lem9} that
\begin{equation}\label{d11}
	\begin{split}
	|\eta_{x_n}(\tilde{x},\tilde{t})|\leq \frac{C}{r} osc_{Q_r^\lambda}\eta\leq C r^{\alpha-1}=C\tilde{x}_n^{\alpha-1},\quad \forall (\tilde{x},\tilde{t})\in K_{1/2}
	\end{split}
\end{equation}
with $r=\tilde{x}_n$.  Notice that the problem  (\ref{a18}) is invariant with respect to the translations in the $x'$ variable. So $\eta_{x_i}(i=1, \cdots, n-1)$ are  also viscosity solutions to   (\ref{a18}).  Thereafter $\eta_{i}:=\partial \eta/\partial x_i $ and $\eta_{i, n}:=\partial^2\eta/(\partial x_i \partial x_n)$.
Then it follows from (\ref{d11}) that 
\begin{equation}\label{d12}
	\begin{split}
		|\eta_{i, n}|\leq Cx_n^{\alpha-1} \quad \text{in} \quad  K_{1/2}.
	\end{split}
\end{equation}
Applying $C^\alpha$ estimates to the solutions $\eta_{i}$ of (\ref{a18}) and by iterating, one gets $\eta\in C^\infty(K_{1/2})$ in the $x'$ variable. Thus we can fix $x'=0$ in  (\ref{a18}). 
For the notation brevity, let
\begin{equation}\label{do12}
	\eta(x,t):=\eta((0, x_n),t)
\end{equation}
 in this step.  We claim that 
\begin{equation}\label{d16}
	\|\nabla \eta\|_{C^{0,\alpha}_{x,t}(K_{1/2}^+)}\leq C,
\end{equation} 
where $C$ is universal constant and the norm $C^{0,\alpha}_{x,t}$ is defined as in (\ref{re2}). In fact it is sufficient to prove pointwise boundary regularity, i.e. on  $x=0$ there is a linear function
$\ell(x)=a\cdot x+b(t)$ such that 
\begin{equation}\label{dv16}
	\|\eta-\ell\|_{L^\infty(Q^{\lambda}_r(0,t))}\leq Cr^{1+\alpha}, \quad\; \forall Q^{\lambda}_r(0,t)\subset K_{1/2}^+,
\end{equation} 
where $Q_r^{\lambda}(0,t)$ is defined as in (\ref{a25}). Together with the interior $C^{1, \alpha}$ estimates of the parabolic equations, we arrive at (\ref{d16}) directly.  \vspace{0.2cm}

Our next aim is to prove (\ref{dv16}). Let
\begin{equation}\label{dr12}
	H(t):=A(\lambda t)PA^T(\lambda t),\quad \;\;v(t):= \nabla\widetilde{\mathcal{G}}(\tilde{a}_n(\lambda t) e_n).
\end{equation}
It is obvious to see that 
$$
\mathcal {L}(A^T(\lambda t) D^2\eta A(\lambda t))=tr(H(t) D^2\eta)=\sum_{i,j=1}^n h^{ij}(t)\eta_{ij}.
$$
By rescaling  $\tilde{\eta}((0, x_n),t):=\eta((0, x_n), \lambda t)$, the problem  (\ref{a18})  is reduced to
\be\label{d13}\le\{\begin{array}{lll}
\tilde{\eta}_t=h^{nn}(t) \tilde{\eta}_{nn}+\sum_{i,j\not=(n,n)}^n h^{ij}(t)\tilde{\eta}_{ij}\quad \text{in}\quad \left(0, \frac{1}{2}\right)\times (-\lambda^{-1}, 0],\\[-0.5ex]
\tilde{\eta}_t=\lambda v_n(t)\tilde{\eta}_n+\lambda \sum_{i=1}^{n-1}v_i(t)\tilde{\eta}_i \quad\quad\;\, \text{on}\quad \{0\}\times (-\lambda^{-1}, 0]
\end{array}\ri.\ee
with the estimate  $(n\Lambda M^2)^{-1}	\leq h^{nn}, v_n\leq \Lambda M^2$, $|h^{nn}_t(t)|\leq \Lambda M^2$,
\begin{equation}\label{d15}
	\Big|\sum_{i,j\not=(n,n)}^n h^{ij}(t)\tilde{\eta}_{ij}\Big|\leq \Lambda M^2 x_n^{\alpha-1}\quad \text{and}\quad 	\Big|\sum_{i=1}^{n-1}v_i(t)\tilde{\eta}_i \Big|\leq \Lambda M^2.
\end{equation} 
 We first show the H$\ddot{\text{o}}$lder estimate of $\nabla\tilde{\eta}$ at $(0,0)$. According to (\ref{af16}), there holds $\lambda\leq \delta$.  Using the  dilation of $\tilde{\eta}$,  one can assume by (\ref{d15}) that 
\begin{equation}\label{df16}
	\Big|\sum_{i,j\not=(n,n)}^n h^{ij}(t)\tilde{\eta}_{ij}\Big|\leq \delta x^{\alpha-1}, \quad \Big|\sum_{i=1}^{n-1}v_i(t)\tilde{\eta}_i \Big|\leq \delta.
\end{equation} 
In what follows, the estimates shall be done in the parabolic cylinders $\tilde{Q}_{\rho}:=(0, \rho)\times (-\rho^2, 0]$.
Suppose
\begin{equation}\label{d19}
	\|\tilde{\eta}-\ell_0\|_{L^{\infty}(\tilde{Q}_{\delta})}\leq C \delta^{1+\alpha},\quad |a_0|\leq 1,
\end{equation} 
where $\ell_0=a_0\cdot x+b_0(t)$ with $b_0'=\lambda v(t)a_0$. Now we want to prove that there is universal  constant $0<\tau<1$ and $\ell_1(x)=a_1x+b_1(t)$ with $|a_1-a_0|\leq C\delta^\alpha$  and $b_1'=\lambda v(t)a_1$ such that
\begin{equation}\label{d20}
	|\tilde{\eta}-\ell_1| \leq C (\tau \delta)^{1+\alpha}\quad \text{in}\quad \tilde{Q}_{\tau\delta}.
\end{equation} 
Then (\ref{dv16}) will be estabilshed by iteration with $\ell_0=0$. To this end, let $\tilde{\eta}$ be a viscosity solution to (\ref{d13}) in $\tilde{Q}_{\delta}$ and
\begin{equation}\label{dff16}
	\gamma(x,t):=\delta^{-(1+\alpha)}(\tilde{\eta}-\ell_0)(\delta x, \delta^2 t).
\end{equation} 
One can easily check that $\gamma$ solves
\be\label{d22}\le\{\begin{array}{lll}
	\gamma_t=h^{nn}(\delta^2 t) \gamma_{nn}+\sum_{i,j\not=(n,n)}^n h^{ij}(\delta^2t)\gamma_{ij}(x,t)-\delta^{1-\alpha} b_0'(\delta^2 t)\quad \text{in}\quad (0, 1)\times (-1, 0],\\[-0.5ex]
	\gamma_t= \lambda \delta \left(v_n(\delta^2 t)\gamma_n+\sum_{i=1}^{n-1}v_i(\delta^2t)\gamma_i\right)\quad\quad\quad\quad\quad\quad\quad\quad\quad \text{on}\quad\{0\}\times (-1, 0].
\end{array}\ri.\ee
It then follows from (\ref{df16}) that
\begin{equation}\label{d21}
	\Big|\sum_{i,j\not=(n,n)}^n h^{ij}(\delta^2t)\gamma_{ij}(x,t)-\delta^{1-\alpha} b_0'(\delta^2 t)\Big|\leq \delta x^{\alpha-1}+C\delta,\quad\;\; \Big|\ \sum_{i=1}^{n-1}v_i(\delta^2 t)\gamma_i \Big|\leq \delta^{1-\alpha}.
\end{equation} 
In view of (\ref{d19}) and (\ref{dff16}), one has $\|\gamma\|_{L^\infty(\tilde{Q}_1)}\leq 1.$ And due to Lemma \ref{lem9}, it yields
$\|\gamma\|_{C^{0,\alpha}_{x,t}(\tilde{Q}_{1/2})}\leq C$. Thus for a sequence $\delta_j\rightarrow 0$, there is a subsequence of functions $\gamma^{(j)}$  such that 
$$\gamma^{(j)}\rightarrow \tilde{\gamma}\in C^{0, \alpha}_{x,t}(\tilde{Q}_{1/2})\quad\text{uniformly}.$$
Thanks  to this,  (\ref{d22}) and (\ref{d21}),  $\tilde{\gamma}$ satisfies 
\begin{equation*}\label{d24}\le\{\begin{array}{lll}
\tilde{\gamma}_t=\tilde{h} \tilde{\gamma}_{nn}\quad \text{in}\quad \Big(0, \frac{1}{2}\Big)\times 	 \Big(-\frac{1}{2}, 0\Big],\\
\tilde{\gamma}_t=0\quad\quad\; \text{on}\quad \{x=0\}\times 	 \Big(-\frac{1}{2}, 0\Big]
\end{array}\ri.\end{equation*} 
in the  viscosity sense,  where $\tilde{h}$ is a constant. Hence $\tilde{\gamma}$ is constant on $\{x=0\}$. By the standard regularity theory, there are constants $\tilde{a}$ and $\tilde{b}$ such that
$$|\tilde{\gamma}-(\tilde{a}x+\tilde{b})|\leq \frac{1}{3}\tau^{1+\alpha}\quad \text{in}\quad  \tilde{Q}_{\tau}\subset K_{1/2}^+$$
for some $0<\tau<1$. This together with (\ref{d19}) and (\ref{dff16})  implies in $\tilde{Q}_{\tau\delta},$
\begin{equation}\label{dgr21}
	|\tilde{\eta}-((a_0+\delta^{\alpha}\tilde{a})x+(b_0(t)+\delta^{1+\alpha}\tilde{b}))|\leq \frac{1}{2}(\tau\delta)^{1+\alpha}.
\end{equation} 
Let $a_1=a_0+\delta^{\alpha}\tilde{a}$ and 
$$
b_1'(t)=\lambda v(t)a_1,\quad 
b_1(0)=b_0(0)+\delta^{1+\alpha}\tilde{b}.$$
It is easy to deduce that
$$|b_1-(b_0(t)+\delta^{1+\alpha}\tilde{b})|\leq \frac{1}{4}(\tau\delta)^{1+\alpha}\quad \text{in}\quad  \tilde{Q}_{\tau\delta}.$$
Combining  (\ref{dgr21}) and the above  leads to 
$$|\tilde{\eta}-(a_1\cdot x+b_1(t))|\leq (\tau\delta)^{1+\alpha}\quad \text{in}\quad  \tilde{Q}_{\tau\delta}.$$
Thus the claim (\ref{d20}) is confirmed. By standard translation arguments,  taking $a_0=0$ and $b_0=0$, we conclude (\ref{dv16}) and hence (\ref{d16}) follows. \vspace{0.2cm}

\emph{Step 2}. In this step, adopting  the notation (\ref{do12}), we shall prove $\eta\in C^1(F_1)$ and 
\begin{equation}\label{dg21}
	\|\nabla \eta\|_{C^{1,\alpha}_{x,t}(K_{1/2}^+)}\leq C.
\end{equation} 
Assume by contradiction that there is a constant $\kappa>0$ and a sequence $t_j\rightarrow0^-$ such that
\begin{equation}\label{d26}
	\eta(0, t_j)-\eta(0,0)>t_j \lambda  v(0)(\eta_x(0,0)-\kappa).
\end{equation} 
We shall compare $\eta$ with the following function in $E_{\kappa}:=[0, c_{\kappa}]\times [t_j, 0]$, 
\begin{equation*}\label{dl27}
	\eta^{sub}(x,t):=\eta(0,0)+x \left(\eta_x(0,0)-\frac{\kappa}{3}\right)+Cx^{1+\alpha}+\lambda  v(0)(\eta_x(0,0)-\kappa) t,
\end{equation*} 
where the constants  $c_{\kappa}$ is small depending on $\kappa$ and $C>0$ is large.  One can check by (\ref{d26}) that $\eta^{sub}$ is a classical subsolution to (\ref{d13}). Then there is a local minimum $(x_0, t_0)\in E_{\kappa}$ such that
\begin{equation}\label{d27}(\eta-\eta^{sub})(x_0, t_0)=\min_{E_\kappa} (\eta-\eta^{sub})(x, t)=0.
\end{equation} 
Clearly,  $x_0\not =0$, since for otherwise, $\eta_x(0, t_0)=0$ and $\eta$ is constant on $\{x=0\}$, which is impossible by Definition \ref{de1}. Thus one has $\{x_0>0\}$. However,  $c_\kappa$ can be chosen small enough such that  $\eta_x>\eta^{sub}_x$. This means that $\eta>\eta^{sub}$ in $(0, c_\kappa]\times [t_j, 0]$,  contradicting to (\ref{d27}). Therefore $\eta\in C^1(F_1).$\vspace{0.2cm}

Next we prove  (\ref{dg21}). For $i\leq  n-1$,  observe  that $\eta_i ((0, x_n), \lambda t)$ solves (\ref{d13}) and $\|D_{x'}^k\eta (x', x_n, t)\|\leq C(k)$ for any $k\geq 1$ in $K_{1/2}$. Thanks to  (\ref{d16}), 
\begin{equation*}\label{dw16}
	\|\nabla \eta_i((0, x_n), t)\|_{C^{0,\alpha}_{x,t}(K_{1/2}^+)}\leq C.
\end{equation*} 
It then follows that 
\begin{equation}\label{dgf21}
	\Big\|\sum_{i,j\not=(n,n)}^n h^{ij}(t)\eta_{ij}\Big\|_{C^{0,\alpha}_{x,t}(K_1)}\leq \Lambda M^2,\quad \Big\|\sum_{i=1}^{n-1} v_i(t)\eta_{i}\Big\|_{C^{\alpha/2}_t(F_1)}\leq \Lambda M^2.
\end{equation} 
Since $\eta\in C^{1, \alpha/2}_{x,t}(F_1)$ and $v_n\in  C^{\alpha/2}_{t}(F_1)$, from (\ref{d13}) and (\ref{dgf21}), 
$$\|\eta\|_{C^{\frac{1+\alpha}{2}}_t(F_1)}\leq C.$$
Applying this and  Schauder estimates to (\ref{d13}), we conclude (\ref{dg21}). Thus  (i) is established.  \vspace{0.2cm}

\emph{Step 3}.  We are going to verify (ii). From Step 1, it is not hard to get
\begin{equation}\label{d30}
	\Big|\eta((0, x_n), t)-\eta(0,t)-\eta_{x_n}(0,t)x_n\Big|\leq Cx_n^{1+\alpha}.
\end{equation} 
Since $\eta$ is smooth with respect to $x'$ variable in $K_{\rho}$,  
\begin{equation}\label{d31}
\eta(x,t)=\eta(0, x_n, t)+\sum_{i<n}\eta_{x_i}(0,0)x_i+O(\rho^{1+\alpha}) \quad \text{for}\quad \rho\leq 1/2. 
\end{equation} 
Combining (\ref{d30}) and (\ref{d31}) gives 
\begin{equation}\label{d32}
|\eta-a(t)x-b(t)|\leq C\rho^{1+\alpha}\quad \text{in}\quad K_{\rho},\quad \forall \rho \leq 1/2,
\end{equation} 
where 
\begin{equation}\label{d33}
a(t)=(\eta_{x_1}(0,0),\ldots, \eta_{x_{n-1}}(0,0),  \eta_{x_n}(0,t)), \quad\; b(t)=\eta(0,t).
\end{equation} 
As a consequence of (\ref{d16}),
\begin{equation}\label{d34}
|a_n(t)-a_n(s)|=|\eta_{x_n}(0,t)-\eta_{x_n}(0,s)|\leq C\lambda^{-\frac{\alpha}{2}}|t-s|^{\frac{\alpha}{2}}.
\end{equation} 
Define a mollifier 
\begin{equation*}\zeta(x)=\le\{\begin{array}{lll}
C \exp\Big({\frac{1}{|x|^2-1}}\Big)\quad \text{if}\quad |x|<1,\\[-0.5ex]
0\quad\quad \quad \quad \quad\;\;\, \text{if}\quad |x|\geq 1,
\end{array}\ri.\end{equation*}
where $C=\Big(\int_{B_1}\exp\Big({\frac{1}{|x|^2-1}}\Big)dx\Big)^{-1}$, and
$$\zeta_{\lambda \rho^2}(x)=\frac{1}{\lambda^n\rho^{2n}}\zeta\left(\frac{x}{\lambda\rho^2}\right).$$
Let the mollification of $a_n$ be given by $\bar{a}_n:=a_n\ast \zeta_{\lambda \rho^2}.$
Let $\bar{a}=(a_1,\ldots, a_{n-1}, \bar{a}_n)$. By scaling back, from (\ref{d34}), one finds   that on  the time interval $\lambda\rho^2$, 
\begin{equation}\label{d38}
	|\bar{a}'(t)|=\lim_{s\rightarrow t}\frac{|\bar{a}(t)-\bar{a}(s)|}{|t-s|}\leq C\lambda^{-\frac{\alpha}{2}}|t-s|^{\frac{\alpha}{2}-1}=C\lambda^{-1}\rho^{\alpha-2},
\end{equation} 
and consequently $|a-\bar{a}|\leq C\rho^\alpha.$ Meanwhile, recall that  $v(t)$ is in (\ref{dr12}). Define $\bar{b}$ satisfying $\bar{b}'(t)=v(t)\bar{a}(t)$ and $\bar{b}(0)=b(0)$. Thanks to (\ref{d33}), one gets
\begin{equation*}\label{d35}
	\begin{split}
		b'=v_n(t)\eta_{x_n}(0,t)+ \sum_{i=1}^{n-1}v_i(t)\eta_{x_i}(0,t)=v(t)\cdot a(t)+\sum_{i=1}^{n-1}v_{i}(t)(\eta_{x_i}(0,t)-\eta_{x_i}(0,0)).
	\end{split}
\end{equation*} 
Then we deduce from above and (\ref{d16}) that
\begin{equation*}\label{d41}
	\begin{split}
		|(\bar{b}'-b')|&\leq|v(t)||(\bar{a}-a)|+ \Big|\sum_{i=1}^{n-1}v_{i}(t)(\eta_i(0,t)-\eta_{x_i}(0,0))\Big|\leq C\rho^\alpha.
	\end{split}
\end{equation*} 
Thus for $t\in(-\rho,0)$, one has $|\bar{b}-b|\leq C\rho^{1+\alpha}$. Then  using (\ref{d32}) we infer 
\begin{equation}\label{dq41}
	\begin{split}
	&\quad |\eta-\bar{a}(t)x-\bar{b}(t)|\\
	&\leq  |\eta-a(t)x-b(t)| +|x||\bar{a}(t)-a(t)|+|\bar{b}(t)-b(t)|\\
	&\leq  C\rho^{1+\alpha}\quad \quad \quad \text{in}\quad K_{\rho}
		\end{split}
\end{equation} 
for each $\rho\leq 1/2$. Furthermore, owing to (\ref{dg21}) and (\ref{d33}), 
$$|a_n(t)-a_n(s)|\leq C\lambda^{-\frac{\alpha+1}{2}}|t-s|^{\frac{\alpha+1}{2}}.$$
By repeating the same arguments as above, one derives $|\bar{a}^{\prime}_n(t)|\leq C \rho^{\alpha-1}\lambda^{-1}$. This together with (\ref{dq41}) confirms (ii).  $\hfill\Box$\\

Finally, we investigate the existence of solutions to  the  following Dirichlet problem. 
\begin{lemma}\label{lem10}
	Let $\psi\in C(\partial_b K_1)$.  Then, there exists a (unique) viscosity solution $\eta$ to 
	\be\label{a32}\le\{\begin{array}{lll}
		\lambda \eta_t=\mathcal {L}(\tilde{A}(\tilde{a}_n(\lambda t) e_n)D^2 \eta )\quad \text{in}\quad K_1,\\[-0.5ex]
		\eta_t=\nabla\widetilde{\mathcal{G}}(\tilde{a}_n(\lambda t) e_n) \nabla \eta\quad\quad \;\; \text{on}\quad F_1,\\[-0.5ex]
		\eta=\psi \quad\quad\quad\quad\quad\quad\quad\quad\;\,\, \text{on}\quad \partial_b K_1.
	\end{array}\ri.\ee	
\end{lemma}
\proof We consider an adaptation of
Perron’s method as argued in \cite{r1}.  Since the proof is identical to that of  \cite[Proposition 5.2]{r1}, we remove the details here, but also refer to \cite{kimm}  for a similar scenario. $\hfill\Box$

\section{Improvement of flatness}	
By the above analysis, this section is devoted   to proving the improvement of flatness lemma, which plays a key role in the proof of Theorem \ref{t1}. We are in the spirit of \cite[Lemma 5.1]{ss1} and \cite[Lemma 5.1]{ss2}.

\begin{lemma}\label{lem1}
	Let $\tilde{u}$ be a viscosity solution to (\ref{a9}). Assume that (\ref{a10})-(\ref{a16}) hold. Then there exists a universal constant  $0<\tau<1$ and a linear function  $l_{\hat{a},\hat{b}}$ with
	\begin{equation}\label{apo12}
		|\tilde{a}(t)-\hat{a}(t)|\leq C \epsilon, \quad  \hat{b}'(t)=\widetilde{\mathcal{G}}(\hat{a}_n(t)e_n),\quad |\hat{a}_n^{\prime}(t)|\leq \frac{\delta\epsilon}{2}(\tau \lambda)^{-2}
\end{equation}
such that
	\begin{equation}\label{a12}
		|\tilde{u}-l_{\hat{a},\hat{b}}(x,t)|\leq \frac{\epsilon}{2} \tau\lambda\quad \text{in}\;\; \overline{K}_{\tau\lambda}.
	\end{equation}
\end{lemma}
\proof We will reason  by contradiction. Suppose that  there is a  sequence $\epsilon_k\rightarrow0$ and viscosity solutions $u^{(k)}$ of (\ref{a9}) with a sequence of operators $\mathcal{F}_k\in \mathscr{S}(\Lambda^{-1}, \Lambda)$, such that
\begin{equation}\label{r1}
	|u^{(k)}-l^{(k)}_{\tilde{a},\tilde{b}}(x,t)|\leq \epsilon_k \lambda\quad \text{in}\;\; \overline{K}_{\lambda},
\end{equation}
where $l^{(k)}_{\tilde{a},\tilde{b}}:=\tilde{a}^{(k)}(t)\cdot x+\tilde{b}^{(k)}(t)$ and $\lambda \leq \min \{\lambda_0, \delta\epsilon_k\}$. Assume in addition that (\ref{a10})-(\ref{asw10}) hold for $u_k$ and  (\ref{aa16}), (\ref{a16}) hold for $\tilde{a}^{(k)}, \tilde{b}^{(k)}$.  However,  $u^{(k)}$  do not satisfy (\ref{a12}).  We decompose the proof into three steps.  And in the proofs, we will omit the superscript of  $\tilde{a}^{(k)}, \tilde{b}^{(k)}$ for the notation simplicity.

\vspace{0.2cm} 

\emph{Step 1} (Compactness for flat sequences)  In this step, we want to prove the following. Under the hypothesises (\ref{r1}),  there is $\mathring{\vartheta}\in C^{0, \alpha}_{x,t}(K^+_{1/2})$  and a subsequence of 
\begin{equation}\label{r3}
	\vartheta^{(k)}\left(x, t \right):=\frac{1}{\epsilon_k \lambda} \left(u^{(k)}(\lambda x, \lambda t)-l^{(k)}_{\tilde{a},\tilde{b}}(\lambda x,\lambda t)\right) \quad \text{in}\quad K_1^+
\end{equation}
such that the following claims hold,\vspace{0.2cm} 

\noindent (i) For every $\delta_0>0$, $\vartheta^{(k)}\rightarrow\mathring{\vartheta}$ uniformly in $K_{1/2}\cap \{x_n\geq \delta_0\}$. \vspace{0.2cm} 

\noindent (ii) The sequence of graphs $\Gamma_k=\{(x, t, \vartheta^{(k)}(x,t)): (x,t)\in K_{1/2}\}$ converges in the Hausdorff distance in $\mathbb{R}^{n+2}$ to the graph $\Gamma=\{(x, t, \mathring{\vartheta}(x,t)): (x,t)\in K_{1/2}^+\}$.\vspace{0.2cm} 

In fact, from Lemma \ref{lem8}, one can see
\begin{equation}\label{rr3}
	|\vartheta^{(k)}(x,t)-\vartheta^{(k)}(y,s)|\leq C (|x-y|^\alpha+\lambda^{-\alpha/2}|t-s|^{\alpha/2})\quad \text{in}\quad Q^{\lambda}_{r}(x,t)\cap K_{1/2}
\end{equation}
with $r\geq C(\delta)\epsilon_k^{1/2} $. Moreover $-1\leq \vartheta_k\leq 1$ in $K_1^+$.  For any given $\delta_0>0$, there is  a sufficiently large $N$ such that when $k\geq N$ one has $C(\delta) \epsilon_{k}\leq \delta_0$. Thus by the Ascoli-Arzel$\grave{\emph{a}}$ theorem, there is a subsequence of $\vartheta^{(k)}$ (which is still denoted by $\vartheta^{(k)}$) such that 
$$\vartheta^{(k)}\rightarrow \mathring{\vartheta}\quad \text{uniformly\; in} \quad K_{1/2}\cap \{x_n\geq \delta_0\},$$ 
where $\mathring{\vartheta}$ is  H$\ddot{\text{o}}$lder continuous over  the set $K_{1/2}\cap \{x_n\geq \delta_0\}$.  Since $C$ and $\alpha$ are independent of $\delta_0$, we can extend $\mathring{\vartheta}$ to  the H$\ddot{\text{o}}$lder continuous function  defined on $K^+_{1/2}$, and thus (i) holds. \vspace{0.2cm}

To prove (ii), we assume 
$X_k=(x_k, t_k, \vartheta^{(k)}(x_k, t_k))\in \Gamma_k.$
Let $k$ be such that $C(\delta) \epsilon_k^{1/2}\leq \delta_0/2$. Let $y_k\in B_{1/2}\cap \{x_n\geq \delta_0\} $ such that $\delta_0/2\leq |x_k-y_k|\leq 2\delta_0$. Then let $Y_k=(y_k, s_k, \vartheta^{(k)}(y_k, s_k))$, and due to (\ref{a25}), one has $|t_k-s_k|\leq 4\lambda\delta_0^2$. This together with (\ref{rr3}) yields
\begin{equation}\label{r6}
		\begin{split}
|X_k-Y_k|^2&\leq |x_k-y_k|^2+\lambda^{-1}|t_k-s_k|+|\vartheta^{(k)}(x_k, t_k)- \vartheta^{(k)}(y_k, s_k)|^2\\
&\leq 8\delta_0^2+ C\delta_0^{2\alpha},
	\end{split}
\end{equation}
where $C$ is universal constant.  In addition, we deduce
\begin{equation}\label{r7}
dist(Y_k, \Gamma)\leq \|\vartheta^{(k)}(y_k, s_k)-\mathring{\vartheta}(y,s)\|_{L^\infty(K_{1/2}\cap \{x_n\geq \delta_0\})}.
	\end{equation}
Combining (\ref{r6}) and (\ref{r7}) leads to
\begin{equation}\label{r8}
	dist(X_k, \Gamma)\leq 2(2\delta_0^2+ C\delta_0^{2\alpha})^{1/2}+\|\vartheta^{(k)}(y_k, s_k)-\mathring{\vartheta}(y,s)\|_{L^\infty(K_{1/2}\cap \{x_n\geq \delta_0\})}.
\end{equation}
On the other hand, assume $\tilde{X}=(x, t, \mathring{\vartheta}(x,t))\in \Gamma$. For each $\delta_0>0$, there exists $y\in B_{\delta_0}(x)\cap (B_{1/2}\cap \{x_n>\delta_0/2\})$. Let $\tilde{Y}=(y, s, \mathring{\vartheta}(y,s))\in \Gamma$ and then  $|t-s|\leq \delta_0^2$. Proceeding as above, as $\mathring{\vartheta}\in C^{0, \alpha}(K^+_{1/2})$,  we estimate 
\begin{equation}\label{r9}
	\begin{split}
		dist(\tilde{X}, \Gamma_k)&\leq |\tilde{X}-\tilde{Y}|+dist(\tilde{Y}, \Gamma_k)\\
		&\leq \left(|x-y|^2+|t-s|+ |\mathring{\vartheta}(x,t)-\mathring{\vartheta}(y,s)|^2\right)^{1/2}+ dist(\tilde{Y}, \Gamma_k)\\
		&\leq C(2\delta_0^2+ \delta_0^{2\alpha})^{1/2}+ \|\vartheta_k(y_k, s_k)-\mathring{\vartheta}(y,s)\|_{L^\infty(K_{1/2}\cap \{x_n\geq \delta_0/2\})}.
	\end{split}
\end{equation}
Finally, in view of the claim (i), (\ref{r8}) and (\ref{r9}), since $\delta_0$ is arbitary, 
$$\lim_{k\rightarrow 0}dist(\Gamma_k, \Gamma)\leq \lim_{k\rightarrow 0}(dist(\tilde{X}, \Gamma_k)+ dist(X_k, \tilde{X}) +dist(X_k, \Gamma))=0. $$ 
This finishes  the proof of the claim (ii) and Step 1. \vspace{0.2cm} 

\emph{Step 2} (Comparision between $\mathring{\vartheta}$ and $\eta$) We first claim that $\mathring{\vartheta}^{sub}:=\mathring{\vartheta}+C\delta(x_n^2-t)$ is a viscosity solution to 
\be\label{df12}\le\{\begin{array}{lll}
	\lambda\mathring{\vartheta}^{sub}_t \leq \mathcal {L}(A^{T}(\lambda t)D^2 \mathring{\vartheta}^{sub} A(\lambda t))\quad \,\text{in}\quad K_1,\\[-0.5ex]
	\mathring{\vartheta}^{sub}_t \leq \nabla\widetilde{\mathcal{G}}(\tilde{a}_n(\lambda t) e_n) \nabla \mathring{\vartheta}^{sub}\quad \quad \quad  \text{on}\quad F_1,
\end{array}\ri.\ee
where the constant $C>0$ is large and universal,   the linear operator $\mathcal {L}$ is given by (\ref{a18})  and  $A(\lambda t):=A(\tilde{a}_n(\lambda t)e_n)$. \vspace{0.2cm} 

To do this, recalling  that  $ \varTheta_n^{(k)}$ is defined in (\ref{a14}), we set 
$$A_k( t):= A(\tilde{a}_n(\lambda t)e_n+\epsilon_k \nabla \vartheta^{(k)}),$$
and denote
$$\widetilde{\mathcal {F}}_k(A_k( t)^T D^2 \vartheta^{(k)} \,A_k(t))=\frac{-\lambda \varTheta_n^{(k)} }{\epsilon_k} \mathcal {F}_k\left(-\frac{\epsilon_k}{\lambda \varTheta_n^{(k)} }  A_k(t)^T D^2 \vartheta^{(k)} \,A_k( t)\right).$$
Since  $M^{-1}\leq \varTheta_n^{(k)}\leq M$ and from (\ref{ff1}), it follows that $\widetilde{\mathcal {F}}_k\in \mathscr{S}(\Lambda^{-1}, \Lambda).$ Thus by extracting a subsequence, together with $\mathcal {F}_k(0)=0$,  one has $\widetilde{\mathcal {F}}_k\rightarrow  \mathcal {L}$ uniformly on compact subsets of matrices as $\epsilon_k \rightarrow 0$. By the assumptions  (\ref{aa16})-(\ref{a16}), 
we conclude  from  Step 1 (i) and Theorem \ref{th4} that $\mathring{\vartheta}^{sub}$ is viscosity solution to $	\lambda \mathring{\vartheta}^{sub}_t \leq \mathcal {L}(A^{T}(\lambda t)D^2 \mathring{\vartheta}^{sub} A(\lambda t))$ in $K_1$. \vspace{0.2cm}

Now we verify that $\mathring{\vartheta}^{sub}$ satisfies the boundary condition of (\ref{df12}) in the viscosity sense. 
Suppose that $P$ is a polynomial touching $\mathring{\vartheta}^{sub}$ from above at $(x_0, t_0)\in F_1$.  Say  $(x_0, t_0)=(0,0)$.  It is clear  that $P$  touches $\mathring{\vartheta}$ above at $(0,0)$. 
We want to prove $P_t(0,0)\leq \nabla\widetilde{\mathcal{G}}(\tilde{a}(0)) \nabla P(0,0)$ for 
\begin{equation}\label{dasw2}
		\lambda P_t > \mathcal {L}(A^{T}(\lambda t)D^2 P A(\lambda t))\quad \text{in}\quad K_r, 
\end{equation}
where $r$ is sufficiently small. Consider a family of polynomials
\begin{equation*}\label{dasw1}
	P^{(\varsigma)}(x, t)=P(x, t)-C\varsigma\delta (\nabla\widetilde{\mathcal{G}}(\tilde{a}(0))+1) (x_n^2+t)+C\varsigma\delta x_n,	\end{equation*}
where $C$ is universal constant and we assume $\nabla\widetilde{\mathcal{G}}(\tilde{a}(0))>0$. If not we replace above by $-\nabla\widetilde{\mathcal{G}}(\tilde{a}(0))$.  Then in a sufficiently small neighborhood of zero,  $P^{(\varsigma)}$  touches  $\mathring{\vartheta}^{sub}$ from above at $(0,0)$. Meanwhile one can deduce 
\begin{equation}\label{dsw1d}
	\begin{split}
		&\quad P^{(\varsigma)}_t(0, 0)-\nabla\widetilde{\mathcal{G}}(\tilde{a}(0)) \nabla P^{(\varsigma)}(0,0)\\
		&=P_t(0, 0)-\nabla\widetilde{\mathcal{G}}(\tilde{a}(0)) \nabla P(0,0) -\tilde{C}\varsigma\delta.
			\end{split}
	\end{equation}
So it suffices to prove that for every $\varsigma>0$, one has
 \begin{equation}\label{dsw2}
	P^{(\varsigma)}_t(0, 0)\leq \nabla\widetilde{\mathcal{G}}(\tilde{a}(0)) \nabla P^{(\varsigma)}(0,0). 
\end{equation}
Fix $\varsigma>0$. Let  $(x_k, t_k)\in K_r^+$ (with $r$ small ) be a sequence of points such that $(x_k, t_k)\rightarrow (0,0)$ as $k\rightarrow \infty$ and  $P^{(\varsigma)}$ touches $\vartheta^{(k)}$  from above at $(x_k, t_k)$. Then we have that $(x_k, t_k)\in F_r$. Suppose on the contrary that $(x_k, t_k)\in K_r$.   One can check that $P^{(\varsigma)}$ is a supersolution to
 (\ref{dasw2}) since $\lambda< \Lambda^{-1}M^{-2}$. Thus  $P^{(\varsigma)}-\vartheta^{(k)}$ is a supersolution to (\ref{a14})  in $K_r$, which makes a contradiction. Consequently our assertion holds. 
Then according to (\ref{r3}),  the polynomial 
$$\hat{P}(x, t)=\epsilon_{k}\lambda P^{(\varsigma)} (x, t)+ l^{(k)}_{\tilde{a},\tilde{b}}(\lambda x,\lambda t)$$
touches $u^{(k)}$ from above at $(x_k, t_k)$. As $u^{(k)}$ is a viscosity solution to (\ref{a9}), we find 
\begin{equation*}\label{dsw3}
	\hat{P}_t(x_k, t_k)-\widetilde{\mathcal{G}}( \nabla \hat{P}(x_k, t_k))\leq 0,
\end{equation*}
and  therefore
$$P^{(\varsigma)}_t(x_k, t_k)\leq \nabla\widetilde{\mathcal{G}}(\tilde{a}_n(\lambda t_k) e_n) \nabla P^{(\varsigma)}(x_k,t_k).$$
Letting $(x_k, t_k)\rightarrow (0,0)$, we conclude (\ref{dsw2}), and then  sending $\varsigma\rightarrow0$ in (\ref{dsw1d}) gives   the desired result.  As a consequence,  the claim  (\ref{df12}) is proven. \vspace{0.2cm}

Note that $\mathring{\vartheta}\in C^{0, \alpha}(K^+_{1/2})$ by Step 1  and $\mathring{\vartheta}^{sub}=\mathring{\vartheta}+C\delta(x_n^2-t)$ solves (\ref{df12}).   Let $\eta=\mathring{\vartheta}^{sub}$ on $\partial_b K_{1/4}$.  By  Lemma \ref{lem9} and  Lemma \ref{lem10}, there exists $\eta\in C^{0, \alpha}_{x,t}(K_{1/4})$ satisfying (\ref{a32}). 
Then it is easy to see
\begin{equation}\label{rg9}
\mathring{\vartheta}^{sub}\leq \eta-C\delta(x_n^2-t-2)=:\tilde{\eta}\quad \text{on}\quad \partial_b K_{1/4}.
\end{equation}
Furthermore, one can check that $\tilde{\eta}$ is a viscosity supersolution to (\ref{a18}). Together with (\ref{rg9}) and Lemma \ref{lm1}, we arrive at $\mathring{\vartheta}-\eta\leq C\delta$ in $K_{1/4}$. Similarly,  one has  $\mathring{\vartheta}^{sup}:=\mathring{\vartheta}-C\delta(x_n^2-t)$ is a viscosity supersolution to (\ref{a18}), i.e.
\begin{equation*}\label{dfs12}\le\{\begin{array}{lll}
	\lambda \mathring{\vartheta}^{sup}_t \geq \mathcal {L}(A^{T}(\lambda t)D^2 \mathring{\vartheta}^{sup} A(\lambda t))\quad \text{in}\quad K_1,\\[-0.5ex]
	\mathring{\vartheta}^{sup}_t \geq \nabla\widetilde{\mathcal{G}}(\tilde{a}_n(\lambda t) e_n) \nabla\mathring{\vartheta}^{sup}\quad \quad \;\;\;\text{on}\quad F_1.
\end{array}\ri.\end{equation*}
Let $\eta=\mathring{\vartheta}^{sup}$ on $\partial_b K_{1/4}$, then 
$$\mathring{\vartheta}^{sup}\geq \eta+C\delta(x_n^2-t-2)=:\hat{\eta}\quad \text{on}\quad \partial_b K_{1/4}.$$
Since $\hat{\eta}$  is a viscosity subsolution to (\ref{a18}), we deduce $\eta-\mathring{\vartheta}\leq C\delta$ in $K_{1/4}$. Thus, 
\begin{equation}\label{gh1}
|\eta-\mathring{\vartheta}|\leq C\delta\quad \text{in}\quad K_{1/4}.
\end{equation}

\emph{Step 3} (Contradiction)  By  (\ref{gh1}), Lemma \ref{lem11} and Step 1, there is a positive integer $N$ such that $k\geq N$, 
$$|\vartheta^{(k)}-l_{\mathring{a}, \mathring{b}}|\leq |\vartheta^{(k)}-\mathring{\vartheta}|+|\mathring{\vartheta}-\eta|+|\eta-l_{\mathring{a}, \mathring{b}}|\leq C\delta+C\rho^{1+\alpha}\quad \text{in}\quad K_{p}. $$ 
Here $\rho\leq 1/4$ and $l_{\mathring{a}, \mathring{b}}$ satisfies (\ref{awq34}). Choose $\rho=\tau, \delta=\tau^{1+\alpha/2}$ and from above
\begin{equation}\label{gh2}
|\vartheta^{(k)}(x,t)-l_{\mathring{a}, \mathring{b}}(x,t)|\leq \frac{1}{4} \tau\quad \text{in}\quad K_{\tau}
\end{equation}
with $|\mathring{a}'_n|\leq \frac{1}{4}\delta\tau^{-2}\lambda^{-1}.$ Owing to (\ref{r3}) and (\ref{gh2}),  
\begin{equation}\label{ghq3}
	\Big|u^{(k)}(x, t)-l^{(k)}_{\hat{a},\check{b}}(x, t)\Big|=\epsilon_{k}\lambda\Big|\vartheta^{(k)}\left(\frac{x}{\lambda}, \frac{t}{\lambda}\right)-l_{\mathring{a}, \mathring{b}}\left(\frac{x}{\lambda}, \frac{t}{\lambda}\right)\Big|\leq  \frac{\epsilon_{k}}{4} \tau\lambda
		\end{equation}
in $K_{\tau\lambda}$, where 
\begin{equation}\label{gh3}
	\begin{split}
		l^{(k)}_{\hat{a},\check{b}}(x, t)& := \hat{a} (t) \cdot x+\check{b}(t)=l^{(k)}_{a,b}(x, t)+\epsilon_{k}\lambda l_{\mathring{a}, \mathring{b}} \left(\frac{x}{\lambda}, \frac{t}{\lambda}\right)\\
		&= \left(\tilde{a}(t)+\epsilon_{k}\mathring{a} \left(\frac{t}{\lambda}\right)\right) \cdot  x+ \left(\tilde{b}(t)+\epsilon_{k}\lambda\mathring{b} \left(\frac{t}{\lambda}\right)\right).
			\end{split}
	\end{equation}
Applying  (\ref{a16}) and (\ref{gh3}), we derive
\begin{equation}\label{ghqo3}
	|\hat{a}'_n(t)|\leq \frac{\epsilon_k \delta}{2(\tau\lambda)^2} \quad \text{and}\quad 	|\hat{a}-\tilde{a}|\leq C\epsilon_{k}.
\end{equation}
Meanwhile, we determine $\hat{b}$ by solving 
$\hat{b}'(t)=\widetilde{\mathcal{G}}(\hat{a}_n(t)e_n)$ and $\hat{b}(0)=\check{b}(0)$. Since 
\begin{equation*}\label{ghs3}
\widetilde{\mathcal{G}}(\hat{a}_n(t)e_n)=\widetilde{\mathcal{G}}\left(\tilde{a}_n(t)e_n+\epsilon_{k}\mathring{a}_n \left(\frac{t}{\lambda}\right)e_n\right) =
\widetilde{\mathcal{G}}(\tilde{a}_n(t)e_n)+\epsilon_{k}\nabla  \widetilde{\mathcal{G}}(\tilde{a}_n(t)e_n) \mathring{a}_n \left(\frac{t}{\lambda}\right)+O(\epsilon_k^2), 
\end{equation*}
when $t\in [-\tau\lambda, 0]$, together with (\ref{gh3}), one deduces 
\begin{equation}\label{gh4}
	\begin{split}
	|\hat{b}(t)-\check{b}(t)|&=\left|\int_{-\tau\lambda}^{0}(\hat{b}'-\check{b}')dt\right|\\
	&\leq \left|\int_{-\tau\lambda}^{0}\left(\hat{b}'-\widetilde{\mathcal{G}}(\tilde{a}_n(t)e_n)-\epsilon_{k}\nabla  \widetilde{\mathcal{G}}(\tilde{a}_n(t)e_n) \mathring{a}_n \left(\frac{t}{\lambda}\right)\right)dt\right|\\
	&= \left|\int_{-\tau\lambda}^{0}\left(\hat{b}'-\widetilde{\mathcal{G}}(\hat{a}_n(t)e_n)+O(\epsilon_{k}^2)\right)dt\right|\leq C\epsilon_{k}^2\tau\lambda.
	\end{split}
\end{equation}
Combining (\ref{ghq3})-(\ref{gh4}), we eventually get 
$$	\Big|u^{(k)}(x, t)-l^{(k)}_{\hat{a},\hat{b}}(x, t)\Big|\leq  \frac{\epsilon_{k}}{2} \tau\lambda$$
with $\hat{a}$ and $\hat{b}$ satisfying (\ref{apo12}). 
This reaches a contradiction to our initial assumption. The lemma is thus proved.  $\hfill\Box$\\

\section{Completion of Proof of Theorem \ref{t1}}	

With the help of Lemma \ref{lem1}, we are going to prove Theorem \ref{t1}  inductively. Then  we shall end this section with a nondegeneracy property of the solution to original  problem (\ref{a1}).   \vspace{0.2cm}

\emph{Completion of proof of Theorem \ref{t1}.} Suppose that  $k\geq 0$ are integers. Denote $\tilde{a}(t)=a_0(t)=(0,\ldots, 0, (a_0)_n(t))$ and $\lambda_k:=\lambda \tau^k$.  Clearly $a_0(t)\in \mathcal{O}_{M/2}$.  Together with (\ref{apo12}), i.e.
\begin{equation}\label{ghd5}
	|\tilde{a}_k(t)-\tilde{a}_{k-1}(t)|\leq C\epsilon_{k}\quad \quad \text{where}\quad \epsilon_{k}:=\epsilon_0 2^{-k},
\end{equation}
it yields $\tilde{a}_k(t)\in \mathcal{O}_M$. Utilizing Lemma \ref{lem1} with $\epsilon=\epsilon_0$,  by induction on $k$,  one gets
$$
|\tilde{u}-l_{\tilde{a}_k, \tilde{b}_k}|\leq \frac{\epsilon_0}{2^k} \tau^k\lambda=\epsilon_{k}\lambda_k.
$$
Setting  $\alpha=-\log\tau^2$ and  choosing $\epsilon_0$ small, we infer 
$\epsilon_{k}\leq (1/2)^k\leq \lambda^{-\alpha}(\lambda\tau^{k})^{\alpha}\leq C(\lambda) \lambda_k^\alpha.$
Therefore, for any $k\geq 0$,  
\begin{equation}\label{gh5}
	|\tilde{u}-l_{\tilde{a}_k, \tilde{b}_k}|\leq C(\lambda) \lambda_k^{1+\alpha}\quad \text{in}\quad K_{\lambda_k}.
\end{equation}
Furthermore,  from (\ref{a17}), one finds for $t, t_0\in [-\lambda_k, 0], $
$$|\nabla \tilde{u}(0, t)-\tilde{a}_k(t_0)|\leq C(\lambda)\lambda_k^{\alpha/2}.$$
Note also by (\ref{ghd5}) that $|\nabla \tilde{u}(0, t)-\tilde{a}_k(t)|\leq C\epsilon_k.$ Hence
$|\tilde{a}_k(t)-\tilde{a}_k(s)|\leq C(\lambda)\lambda_k^{\alpha/2}$ for $t, s\in [-\lambda_k, 0].$ This together with (\ref{gh5}) and the facts that 
$b'(t)=\widetilde{\mathcal{G}}(a_n(t)e_n)$
leads to
$$	|\tilde{u}-(\tilde{a}_k(0)\cdot x+\tilde{b}'_k(0)t+\tilde{b}_k(0))|\leq C(\lambda) \lambda_k^{1+\alpha/2}\quad \text{in}\quad K_{\lambda_k}.$$

Now we go back to the statements of Theorem \ref{t1}.  From the  $\bar{\epsilon}$-flat  assumption, there holds
$$|u(x, t)-(a_n(t)x_n-b(t))^+|\leq \bar{\epsilon}\lambda.$$
Set $\epsilon=\bar{\epsilon}/\tau$ and $\bar{\lambda}=\tau \lambda$, where $\tau\leq(\delta \bar{\epsilon}/\lambda)^{1/2}$. So  $\bar{\lambda}\leq \min \{\lambda_0, \delta \epsilon\}$. By working in the cylinders $K_{\bar{\lambda}}$, we may repeat the previous argument to get the desired result. That is,  there are universal constants $\alpha\in(0,1)$ and $C>0$ such that for all $(x',x_n, t), (y',y_n, s)\in \mathcal{Q}_{\lambda/2}$, the free boundary
$$\partial\Omega^+_u=\{(x', x_n, t) \in \mathcal{Q}_{\lambda/2}:  x_n=\tilde{u}(x',t)\}$$
with the estimate
$$|\tilde{u}(x', t)-\tilde{u}(y',s)-\nabla \tilde{u} (y',s)(x'-y')|\leq C (|x'-y'|+\sqrt{|t-s|})^{1+\alpha}.$$ 
Hence $\partial\Omega^+_u$ is $C^{1,\alpha}$ in the $x_n$ direction. \vspace{0.2cm}

Next applying the Schauder estimates (\cite[Corollary 14.9]{gl} and \cite[Proposition 5.3]{goffi}) to (\ref{a9}), we reach $C^{\infty}$ regularity of $\tilde{u}$. Therefore, Theorem \ref{t1} is  concluded.  $\hfill\Box$\\

\emph{Nondegeneracy property of $u$.}   For some $x_0\in \bar{B}_{3\lambda/4}$, we assert that  $u(x_0,t)\geq c_1\lambda$ with $t\in [-c_2\lambda, 0]$, where $c_1, c_2>0$ are some constants. Indeed, assume $x_0=\frac{3}{4}\lambda e_n$. Since $(0,0)\in \partial\Omega^+_u$, from (\ref{a2}) and (\ref{ar2}), we have $b(0)<\bar{\epsilon}\lambda$.  Moreover,  the conditions (\ref{addf1})  and (\ref{ar2})  imply $0<\tilde{M}^{-1}\leq \mathcal{G} (a_n)\leq \tilde{M}$. Thus for $t>-c_2\lambda=-\frac{\lambda}{3M \tilde{M}}$, there holds $b(t)<\bar{\epsilon}\lambda-t\tilde{M}<\frac{\lambda }{2M}$ for small $\bar{\epsilon}$. 
Note also that  $u$ is $\bar{\epsilon}$-flat.  	It then follows that
		\begin{equation*}
			\begin{split}
				u\geq \left(\frac{3\lambda}{4M}-b(t)-\bar{\epsilon}\lambda\right)^+\geq \left(\frac{5\lambda}{8M}-b(t)\right)^+\geq \frac{\lambda}{8M},
			\end{split}
		\end{equation*}
as desired. $\hfill\Box$

\vspace{0.8cm}
\noindent \textbf{Acknowledgements.} The author  would like to thank  Prof. Hui Yu for many helpful discussions and stimulating conversations on the direction of  this topic. 
\vspace{0.2cm}
%\bigskip

\end{document}